\documentclass{article}
\usepackage{amsmath}
\usepackage{amssymb}

\input{tcilatex}
\begin{document}

\begin{center}
\textbf{Abstract elliptic and parabolic equat\i ons in Morrey spaces and
applications} 

\textbf{$M.A.$ $Ragusa^{1,2,\ast }$}, \textbf{$V.B.$ $Shakhmurov^{3}$}

1. Department of Mathematics, University of Catania, Viale Andrea Doria n.6,
Catania 95128, Italy \\[0pt]
2. RUDN University, Miklukho-Maklaya Str. 6, Moscow, 117198, Russia,

E-Mail: maragusa@dmi.unict.it

3. Department of Mechanical Engineering, Okan University, Akfirat, Tuzla
34959 Istanbul, Turkey,

E-mail: veli.sahmurov@okan.edu.tr

\textbf{\ Abstract}
\end{center}

We presents the study the separability properties for differential-operator
equations in Morrey spaces. We prove that the corresponding differential
operator is a generator of analytic semigroup in vector-valued Morrey
spaces. Moreover, \ maximal regularity properties of corresponding parabolic
equation \i s obtained. In applications, the maximal regularity properties
of \ Wentzell-Robin type problem for elliptic equations and mixed value
problem for degenerate parabolic equations in Morrey spaces are derived.


\begin{center}
\bigskip\ \ \textbf{AMS:34G10, 35J25, 35J70}

\textbf{1. Introduction,\ notations and background}
\end{center}

The goal of the present paper is to study the boundary value problems (BVPs)
for differential-operator equations (DOEs) with VMO top-order coefficients

\begin{equation}
\Sigma_{k=1}^{n}a_{k}\left( x_{k}\right) \frac{\partial ^{2}u}{\partial
x_{k}^{2}}+A\left( x\right) u+\Sigma_{k=1}^{n}A_{k}\left( x\right) \frac{%
\partial u}{\partial x_{k}}+\mu u=f\left( x\right) ,\text{ }x\in \mathbb{R}%
_{+}^{n},  \tag{1.1}
\end{equation}

\begin{equation*}
L_{1}u=\sum\limits_{i=0}^{\sigma }\alpha _{i}u^{\left( i\right) }\left(
x^{\prime },0\right) =0\text{, }x^{\prime }=\left(
x_{1},x_{2},...,x_{n-1}\right) ,\text{ }\sigma \in \left\{ 0,1\right\} ,
\end{equation*}%
\ \ where $a_{k}$ are complex-valued functions, $\alpha _{i}$ are complex
numbers, $A=A\left( x\right) ,$ $A_{k}=A_{k}\left( x\right) $ are linear
operators in a Banach space $E$ and $\mu $ is a complex parameter$.$

It is known that many classes of PDEs, pseudo DEs and integro DEs can be
expressed in the form of DOEs. Therefore, many researchers (see e.g. $\left[ 
\text{1}\right] $, $\left[ \text{5-6}\right] $, $\left[ \text{15}\right] ,$ $%
\left[ \text{24-27}\right] ,\left[ \text{30}\right] $) investigated similar
classes of PDEs under a single DOE. Moreover, the maximal regularity
properties of DOEs with continuous coefficients were studied e.g. in $\left[
5\right] $, $\left[ 24-27\right] $ and $\left[ 30\right] .$ Moreover, the
regularity propert\i es of elliptic and parabolic equations in Morrey spaces
were defined e.g. in $\left[ 7-9\right] $, $\left[ 13-14\right] ,$ $\left[
19-22\right] .$

Here the equation with operator coefficients is considered in abstract
Morrey spaces. We shall prove separability of the problem $(1)$, i.e. we
show that for each $f\in L^{p,\lambda }\left( \mathbb{R}_{+}^{n};E\right) $
there exists a unique strong solution $u$ of the problem $\left( 1.1\right) $
and a positive constant $C$ depending only on $p,$ $E$ and $A$ such that 
\begin{equation*}
\sum\limits_{i=0}^{2}\Sigma_{k=1}^{n}\left\vert \mu \right\vert ^{1-\frac{i}{%
2}}\left\Vert \frac{\partial ^{i}u}{\partial x_{k}}\right\Vert
_{L^{p,\lambda }\left( \mathbb{R}_{+}^{n};E\right) }+\left\Vert
Au\right\Vert _{L^{p,\lambda }\left( \mathbb{R}_{+}^{n};E\right) }\leq
C\left\Vert f\right\Vert _{L^{p,\lambda }\left( \mathbb{R}_{+}^{n};E\right)
}.
\end{equation*}

Then, the well-posedness of mixed problem for the parabolic equation 
\begin{equation}
\frac{\partial u}{\partial t}+\Sigma _{k=1}^{n}a_{k}\left( x_{k}\right) 
\frac{\partial ^{2}u}{\partial x_{k}^{2}}+A\left( x\right) u=f\left(
x,t\right) \text{, }x\in \mathbb{R}_{+}^{n}\text{, }t\in \left( 0,T\right) ,
\tag{1.2}
\end{equation}%
\begin{equation*}
\sum\limits_{i=0}^{\sigma }\alpha _{i}u_{x_{n}}^{\left( i\right) }\left(
x^{\prime },0,x\right) =0,x^{\prime }=\left( x_{1},x_{2},...,x_{n-1}\right) 
\text{, }\sigma \in \left\{ 0,1\right\} ,
\end{equation*}%
\begin{equation*}
u\left( x,0\right) =0
\end{equation*}%
is derived in $E$-valued mixed Morrey space $L^{p,\lambda }\left( \mathbb{R}%
_{+}^{n}\times \left( 0,T\right) ;E\right) $.\ \ \ \ \ \ \ \ \ \ 

Note that, the principal part of corresponding differential operator is non
selfadjoint. Nevertheless, the sharp uniform coercive estimate for\ the
resolvent, Fredholmness are established. In application, we study
Wentzell-Robin type problem for elliptic equations and mixed value problems
for degenerate parabolic equations with $VMO$ coefficients in Morrey spaces.

Since $\left( 1\right) $\ involves unbounded operators, it is not easy to
get representation for Green function and estimate of solutions$.$ Therefore
we use the modern harmonic analysis elements e.g. the Hilbert operator and
the commutator estimates in $E$-valued $L_{p}$ spaces, embedding theorems of
Sobolev-Lions spaces and some semigroups estimates to overcome these
difficulties. Moreover we also use our previous results on equations with
continuous leading coefficients and the perturbation theory of linear
operators to obtain our main assertions.

Since the Hilbert space $E$ is arbitrary and $A$ is a possible linear
operator, by choosing $E$ and $A$ we can obtain numerous classes of elliptic
and parabolic type equations and its systems which occur in the different
processes. \ In application, we put $E=L^{p_{1}}\left( 0,1\right) $ and $A$
to be differential operator with generalized Wentzell-Robin boundary
condition defined by 
\begin{equation*}
D\left( A\right) =\left\{ u\in W_{p_{1}}^{2}\left( 0,1\right) ,\text{ }%
B_{j}u=Au\left( j\right) ,\text{ }j=0,1\right\} ,\text{ }
\end{equation*}%
\begin{equation*}
\text{ }Au=au^{\left( 2\right) }+bu^{\left( 1\right) },
\end{equation*}%
in $\left( 1.1\right) -\left( 1.2\right) $, where $a,$ $b$ are
complex-valued functions. Then, we obtain the following Wentzell-Robin type
mixed problem for elliptic equations

\begin{equation}
\Sigma _{k=1}^{n}a_{k}\left( x_{k}\right) \frac{\partial ^{2}u}{\partial
x_{k}^{2}}+A\left( x\right) u+a\frac{\partial ^{2}u}{\partial y^{2}}+b\frac{%
\partial u}{\partial y}=f\left( x,y\right) ,\text{ }x\in \mathbb{R}_{+}^{n},
\tag{1.3}
\end{equation}

\begin{equation}
\sum\limits_{i=0}^{\sigma }\alpha _{i}\frac{\partial ^{i}u}{\partial
x_{n}^{i}}\left( x^{\prime },0\right) =0\text{, }\sigma \in \left\{
0,1\right\} \text{, }u=u\left( x,y\right) ,  \tag{1.4}
\end{equation}%
\begin{equation}
Au\left( x,j\right) =0\text{, }j=0,1,\text{ }u\left( x,0\right) =a\left(
x\right) ,\text{ }x^{\prime }\in \mathbb{R}^{n-1}\text{, }y\in \left(
0,1\right) .  \tag{1.5}
\end{equation}%
\ \ \ 

Note that, the regularity properties of Wentzell-Robin type boundary value
problems (BVP) for elliptic equations were studied e.g. in $\left[ \text{10,
11}\right] $ and the references therein. Here, 
\begin{equation*}
\tilde{\Omega}=\mathbb{R}_{+}^{n}\times \left( 0,1\right) \text{, }\bar{%
\Omega}=\mathbb{R}_{+}^{n}\times \left( 0,b\right) \text{, }\mathbf{p=}%
\left( p_{1},p\right) \text{, }p_{1},p\in \left( 1,\infty \right)
\end{equation*}%
$L^{\mathbf{p,}\lambda }\left( \tilde{\Omega}\right) $ denotes the space of
all $\mathbf{p}$-summable Morrey space with mixed norm i.e., the space of
all measurable functions $f$ defined on $\tilde{\Omega}$, for which 
\begin{equation*}
\left\Vert f\right\Vert _{L^{\mathbf{p,}\lambda }\left( \tilde{\Omega}%
\right) }=\left( \int\limits_{R_{+}^{n}}\sup_{x_{0}\mathbb{R}_{+}^{n}\in
,r}r^{-\lambda }\left( \int\limits_{0}^{1}\left\vert f\left( x,y\right)
\right\vert ^{p_{1}}dy\right) ^{\frac{p}{p_{1}}}dx\right) ^{\frac{1}{p}%
}<\infty \text{, }0<\lambda <n.
\end{equation*}

By using the general abstract result above, the existence, uniqueness and
maximal $L^{\mathbf{p,}\lambda }\left( \tilde{\Omega}\right) $ regularity
properties of the problem $\left( 1.3\right) -\left( 1.5\right) $ is
obtained.\ 

Moreover, let we choose $E=L^{^{p_{1}}}\left( 0,b\right) $ and $A$ to be
degenerated differential operator in $L^{p_{1}}\left( 0,b\right) $ defined
by 
\begin{equation*}
D\left( A\right) =\left\{ u\in W_{\gamma }^{\left[ 2\right] ,p_{1}}\left(
0,1\right) ,\right. \left. \dsum\limits_{i=0}^{\nu _{k}}\alpha _{ki}u^{\left[
i\right] }\left( 0\right) +\beta _{ki}u^{\left[ i\right] }\left( b\right) =0,%
\text{ }k=1,2\right\} ,\text{ }
\end{equation*}%
\begin{equation}
\text{ }A\left( x\right) u=b_{1}\left( x,y\right) u^{\left[ 2\right]
}+b_{2}\left( x,y\right) u^{\left[ 1\right] }\text{, }x\in \mathbb{R}%
_{+}^{n},\text{ }y\in \left( 0,b\right) ,\text{ }\nu _{k}\in \left\{
0,1\right\} ,  \tag{1.6}
\end{equation}%
\ \ \ where $u^{\left[ i\right] }=\left( y^{\gamma }\frac{d}{dy}\right)
^{\gamma }u$ for $0\leq \gamma <\frac{1}{2},$ $b_{1}=b_{1}\left( x,y\right) $
is a cont\i nous, $b_{2}=b_{2}\left( x,y\right) $ is a bounded functon on $%
y\in $ $\left[ 0,1\right] $ for a.e. $x\in \mathbb{R}_{+}^{n},$ $\alpha
_{ki} $, $\beta _{ki}$ are complex numbers and $W_{\gamma }^{\left[ 2\right]
,p_{1}}\left( 0,b\right) $ is a weighted Sobolev spase defined by 
\begin{equation*}
W_{\gamma }^{\left[ 2\right] ,p_{1}}\left( 0,b\right) =\left\{ {}\right.
u:u\in L^{p_{1}}\left( 0,b\right) ,\text{ }u^{\left[ 2\right] }\in
L^{p_{1}}\left( 0,b\right) ,\text{ }
\end{equation*}%
\begin{equation*}
\left\Vert u\right\Vert _{W_{\gamma }^{\left[ 2\right] ,p_{1}}}=\left\Vert
u\right\Vert _{L^{p_{1}}}+\left\Vert u^{\left[ 2\right] }\right\Vert
_{L^{p_{1}}}<\infty .
\end{equation*}%
Then, from $\left( 1.2\right) $ we get the following mixed problem for
degenerate parabolic equation

\begin{equation}
\frac{\partial u}{\partial t}+\Sigma _{k=1}^{n}a_{k}\left( x_{k}\right) 
\frac{\partial ^{2}u}{\partial x_{k}^{2}}+\left( b_{1}\frac{\partial ^{\left[
2\right] }u}{\partial y^{2}}+b_{2}\frac{\partial ^{\left[ 1\right] }u}{%
\partial y}\right) =f\left( x,y,t\right) ,\text{ }  \tag{1.7}
\end{equation}%
\begin{equation*}
\text{ }y\in \left( 0,b\right) ,\text{ }t\in \left( 0,T\right) ,\text{ }%
u=u\left( x,y,t\right) ,
\end{equation*}%
\begin{equation}
\sum\limits_{i=0}^{\sigma }\alpha _{i}\frac{\partial ^{i}u}{\partial
x_{n}^{i}}\left( x^{\prime },0,y,t\right) =0\text{, }\sigma \in \left\{
0,1\right\} \text{, }u=u\left( x,y,t\right) ,  \tag{1.8}
\end{equation}%
\ \ \ 

\begin{equation*}
L_{1}u=\dsum\limits_{i=0}^{\nu _{k}}\alpha _{ki}u^{\left[ i\right] }\left(
x,0,t\right) +\beta _{ki}u^{\left[ i\right] }\left( x,b,t\right) =0,\text{ }%
k=1,2,
\end{equation*}

\begin{equation}
u\left( x,y,0\right) =a\left( x,y\right) .  \tag{1.9}
\end{equation}

\bigskip The existence, uniqueness and coercive $L^{\mathbf{p}}\left( \bar{%
\Omega}\right) $ estimates for solution of the problem $\left( 1.7\right)
-\left( 1.9\right) $ is derived.

\bigskip\ \ \ Let $E$ be a Banach space and $\gamma =\gamma \left( x\right)
, $ $x=\left( x_{1},x_{2},...,x_{n}\right) $ be a positive measurable
function on a on the measurable subset $\Omega \subset R^{n}.$ Let $%
L_{p,\gamma }\left( \Omega ;E\right) $ denote the space of strongly
measurable $E-$valued functions that are defined on $\Omega $ with the norm

\begin{equation*}
\left\Vert f\right\Vert _{L_{\gamma }^{p}\left( E\right) }=\left\Vert
f\right\Vert _{L_{\gamma }^{p}\left( \Omega ;E\right) }=\left(
\int\limits_{\Omega }\left\Vert f\left( x\right) \right\Vert _{E}^{p}\gamma
\left( x\right) dx\right) ^{\frac{1}{p}},1\leq p<\infty .
\end{equation*}

For $\gamma \left( x\right) \equiv 1,$ the space $L_{\gamma }^{p}\left(
\Omega ;E\right) $ will be denoted by $L^{p}\left( E\right) =L^{p}\left(
\Omega ;E\right) .$ Let $B(x,r)$ denote the ball of center $x$ and radius $r$%
. Let $1<p<\infty $ and $0<\lambda <n.$ \textbf{\ }By $L^{p,\lambda }\left( 
\mathbb{R}^{n};E\right) $ we will denote $E$-valued Morrey space that is the
set of all functions $f\in L_{loc}(\mathbb{R}^{n};E)$ for which 
\begin{equation*}
||f||_{L^{p,\lambda }(\mathbb{R}^{n};E)}^{p}=\sup_{x\in \mathbb{R}%
^{n},r}r^{-\lambda }\int\limits_{B(x,r)}\left\Vert f\left( y\right)
\right\Vert _{E}^{p}dy<\infty .
\end{equation*}

$C(\Omega ,E)$ and $C^{m}(\Omega ;E)$ will denote the spaces of $E$-valued
bounded uniformly strongly continuous and $m$-times continuously
differentiable functions on $\Omega ,$ respectively. $C_{0}^{m}(\Omega ;E)$%
-denotes the class of $E$-valued functions from $C^{m}(\Omega ;E)$ with
commpact supports on $\Omega .$

Let $E_{0}$ and $E$ be two Banach spaces and $E_{0}$ continuously and
densely embedded into $E.$ Let $m$ be a natural number and $\Omega \subset 
\mathbb{R}^{n}.$Let$\ \Omega $ be a domain on $R^{n}$ and. Here, $%
W^{m,p}\left( \Omega ;E\right) $ denotes the abstract Sobolev space of
functions $u\in L^{p}\left( \Omega ;E\right) $ which have generalized
derivatives $D^{\alpha }u\in L_{p}\left( \Omega ;E\right) $ for $D^{\alpha
}=D_{x_{1}}^{\alpha _{1}}D_{x_{2}}^{\alpha _{2}},...,D_{x_{n}}^{\alpha _{n}}$
and $\alpha =\left( \alpha _{1},\alpha _{2},...,\alpha _{n}\right) $\ with
norm 
\begin{equation*}
\left\Vert u\right\Vert _{W^{m,p}\left( \Omega ;E\right)
}=\sum\limits_{\left\vert \alpha \right\vert \leq m}\left\Vert D^{\alpha
}u\right\Vert _{L^{p}\left( \Omega ;E\right) }<\infty .
\end{equation*}

\ \ $W^{m,p}\left( \Omega ;E_{0},E\right) $ denotes the Sobolev-Lions type
space of functions $u\in L^{p}\left( \Omega ;E_{0}\right) $ which have
generalized derivatives $D^{\alpha }u\in L^{p}\left( \Omega ;E\right) $ with
norm 
\begin{equation*}
\left\Vert u\right\Vert _{W^{m,p}\left( \Omega ;E_{0},E\right) }=\left\Vert
u\right\Vert _{L^{p}\left( \Omega ;E_{0}\right) }+\sum\limits_{\left\vert
\alpha \right\vert \leq m}\left\Vert D^{\alpha }u\right\Vert _{L^{p}\left(
\Omega ;E\right) }<\infty .
\end{equation*}

It is clear to see that 
\begin{equation*}
W^{m,p}\left( \Omega ;E_{0},E\right) =L^{p}\left( \Omega ;E_{0}\right) \cap
W^{m,p}\left( \Omega ;E\right) .
\end{equation*}

The Banach space\ $E$\ is called a UMD-space and written as $E\in $ UMD if
only if the Hilbert operator 
\begin{equation*}
\left( Hf\right) \left( x\right) =\lim\limits_{\varepsilon \rightarrow
0}\int\limits_{\left\vert x-y\right\vert >\varepsilon }\frac{f\left(
y\right) }{x-y}dy
\end{equation*}%
is bounded in the space $L_{p}\left( R,E\right) ,$ $p\in \left( 1,\infty
\right) $ (see e.g. $\left[ 3\right] $). UMD spaces include $L_{p}$, $l_{p}$
spaces, Lorentz spaces $L_{pq},$ $p,$ $q\in \left( 1,\infty \right) $ and
Morrey spaces (see e.g. $\left[ 20\right] $).

Here $L^{p,\lambda }\left( \mathbb{\Omega };E\right) $ wii denote $E$-valued
Morrey space that is the set of all functions $f\in L_{loc}(\mathbb{\Omega }%
;E)$ for which 
\begin{equation*}
||f||_{L^{p,\lambda }(\mathbb{\Omega };E)}^{p}=\sup_{x_{0}\in \Omega
,r}r^{-\lambda }\int\limits_{B(x,r)\cap \Omega }\left\Vert f\left( y\right)
\right\Vert _{E}^{p}dy<\infty .
\end{equation*}

Let \ 
\begin{equation*}
S\left( \varphi \right) =\left\{ \lambda \in \mathbb{C}\text{, }\left\vert
\arg \lambda \right\vert \leq \varphi \right\} \cup \left\{ 0\right\} ,0\leq
\varphi <\pi .
\end{equation*}

A linear operator\ $A$ is said to be a $\varphi $-positive in a Banach\
space $E$ with bound $M>0$ if $D\left( A\right) $ is dense on $E$ and $%
\left\Vert \left( A+\lambda I\right) ^{-1}\right\Vert _{L\left( E\right)
}\leq M\left( 1+\left\vert \lambda \right\vert \right) ^{-1}$ $\ $with $%
\lambda \in S\left( \varphi \right) ,\varphi \in \left( 0,\pi \right] $, $I$
is an identity operator in $E$ and $L\left( E\right) $ is the space of
bounded linear operators in $E.$ Sometimes instead of $A+\lambda I$\ will be
written $A+\lambda $ and it denoted by $A_{\lambda }.$ It is known $\left[ 
\text{29, \S 1.15.1}\right] $ that there exists fractional powers\ $%
A^{\theta }$of the sectorial operator $A.$ Let $E\left( A^{\theta }\right) $
denote the space $D\left( A^{\theta }\right) $ with graphical norm 
\begin{equation*}
\left\Vert u\right\Vert _{E\left( A^{\theta }\right) }=\left( \left\Vert
u\right\Vert ^{p}+\left\Vert A^{\theta }u\right\Vert ^{p}\right) ^{\frac{1}{p%
}},1\leq p<\infty ,\text{ }-\infty <\theta <\infty .
\end{equation*}

\ Let $E_{1}$ and $E_{2}$ be two Banach spaces. By $\left(
E_{1},E_{2}\right) _{\theta ,p}$, $0<\theta <1,1\leq p\leq \infty $ will be
denoted the interpolation spaces obtained from $\left\{ E_{1},E_{2}\right\} $
by the $K$-method $\left[ \text{29, \S 1.3.1}\right] $.\ 

A set $W\subset B\left( E\right) $ is called $R$-bounded (see e.g. $\left[ 
\text{5, \S\ 3.1}\right] $ or $\left[ \text{28}\right] $ ) if there is a
positive constant $C$ such that for all $T_{1},T_{2},...,T_{m}\in W$ and $%
u_{1,}u_{2},...,u_{m}\in E_{1},$ $m\in \mathbb{N}$ 
\begin{equation*}
\int\limits_{0}^{1}\left\Vert \sum\limits_{j=1}^{m}r_{j}\left( y\right)
T_{j}u_{j}\right\Vert _{E}dy\leq C\int\limits_{0}^{1}\left\Vert
\sum\limits_{j=1}^{m}r_{j}\left( y\right) u_{j}\right\Vert _{E}dy,
\end{equation*}

where $\left\{ r_{j}\right\} $ is a sequence of independent symmetric $%
\left\{ -1,1\right\} $-valued random variables on $\left[ 0,1\right] $.

Let $S\left( \mathbb{R}^{n};E\right) $ denote the Schwartz class, i.e. the
space of all $E$-valued rapidly decreasing smooth functions on $R.$ Let $F$
denote the Fourier transformation. A function $\Psi \in L^{\infty }\left( 
\mathbb{R}^{n};B\left( E\right) \right) $ is called a Fourier multiplier
from $L^{p}\left( \mathbb{R}^{n};E\right) $\ to $L^{p}\left( \mathbb{R}%
^{n};E\right) $ if the map $u\rightarrow \Lambda _{\Psi }u=F^{-1}\Psi \left(
\xi \right) Fu,$ $u\in S\left( \mathbb{R}^{n};E\right) $ is well defined and
extends to a bounded linear operator 
\begin{equation*}
\Lambda _{\Psi }:\ L^{p}\left( \mathbb{R}^{n};E\right) \rightarrow \
L^{p}\left( \mathbb{R}^{n};E\right) .
\end{equation*}%
The set of all multipliers from\ $L^{p}\left( \mathbb{R}^{n};E\right) $ to\ $%
L^{p}\left( \mathbb{R}^{n};E\right) $\ will be denoted by $M_{p}^{p}\left(
E\right) .$

\textbf{Definition 1.1}$.$ The $\varphi $-positive operator $A$ is said to
be an $R$-positive in a Banach space $E$ if there exists $\varphi \in \left[
0\right. \left. \pi \right) $ such that the set 
\begin{equation*}
\left\{ A\left( A+\lambda \right) ^{-1}:\lambda \in S\left( \varphi \right)
\right\}
\end{equation*}%
is $R$-bounded.

A linear operator $A\left( x\right) $ is said to be positive in $E$
uniformly in $x$ if $D\left( A\left( x\right) \right) $ is independent of$\
x $, $D\left( A\left( x\right) \right) $ is dense in $E$ and $\left\Vert
\left( A\left( x\right) +\lambda \right) ^{-1}\right\Vert \leq M\left(
1+\left\vert \lambda \right\vert \right) ^{-1}$ for all $\lambda \in S\left(
\varphi \right) $, $\varphi \in \left[ 0\right. \left. \pi \right) .$

\vspace{3mm} Let $\sigma _{\infty }\left( E_{1},E_{2}\right) $ denote the
space of all compact operators from $E_{1}$ to $E_{2}.$ For $E_{1}=E_{2}=E$
it is denoted by $\sigma _{\infty }\left( E\right) .$

Let $m$ be a natural number and $\Omega \subset \mathbb{R}^{n}.$\ Let $%
W^{m,p}\left( \Omega ;E_{0},E\right) $ denote a space of all $\ $functions $%
u\in L^{p}\left( \Omega ;E_{0}\right) $ possess the generalized derivatives $%
D^{\alpha }u$ such that $D^{\alpha }u\in L^{p}\left( \Omega ;E\right) $ with
the norm 
\begin{equation*}
\left\Vert u\right\Vert _{W^{m,p}\left( \Omega ;E_{0},E\right) }=\left\Vert
u\right\Vert _{L^{p}\left( \Omega ;E_{0}\right) }+\sum\limits_{\left\vert
\alpha \right\vert \leq m}\left\Vert D^{\alpha }u\right\Vert _{L^{p}\left(
\Omega ;E\right) }<\infty .
\end{equation*}%
\ This called Sobolev-Lions type space. For $E_{0}=$ $E$ the space $%
W^{m,p}\left( \Omega ;E_{0},E\right) $ will denoted by $W^{m,p}\left( \Omega
;E\right) .$ It is clear to see that 
\begin{equation*}
W^{m,p}\left( \Omega ;E_{0},E\right) =W^{m,p}\left( \Omega ;E\right) \cap
L^{p}\left( \Omega ;E_{0}\right) .
\end{equation*}

Now, consider Sobolev- Morrey-Lions type space $W^{m,p,\lambda }\left(
\Omega ;E_{0},E\right) $ that is a space of all $\ $functions $u\in
L^{p,\gamma }\left( \Omega ;E_{0}\right) $ possess the generalized
derivatives $D^{\alpha }u$ such that $D^{\alpha }u\in L^{p,\lambda }\left(
\Omega ;E\right) $ with the norm 
\begin{equation*}
\left\Vert u\right\Vert _{W^{m,p,\lambda }\left( \Omega ;E_{0},E\right)
}=\left\Vert u\right\Vert _{L^{p,\lambda }\left( \Omega ;E_{0}\right)
}+\sum\limits_{\left\vert \alpha \right\vert \leq m}\left\Vert D^{\alpha
}u\right\Vert _{L^{p,\lambda }\left( \Omega ;E\right) }<\infty .
\end{equation*}%
\ For $E_{0}=$ $E$ the space $W^{m,p,\gamma }\left( \Omega ;E_{0},E\right) $
will denoted by $W^{m,p,\lambda }\left( \Omega ;E\right) .$ It is clear to
see that 
\begin{equation*}
W^{m,p,\gamma }\left( \Omega ;E_{0},E\right) =W^{m,p,\lambda }\left( \Omega
;E\right) \cap L^{p,\lambda }\left( \Omega ;E_{0}\right) .
\end{equation*}

Function $u\in W^{2,p,\lambda }\left( 0,1;E\left( A\right) ,E\right) $
satisfying the equation $\left( 1.1\right) $ a.e. on $\left( 0,1\right) $ is
said to be solution of the problem $\left( 1.1\right) .$

From $\left[ 24\right] $ we obtain:

\textbf{Theorem A}$_{1}.$ Let the following conditions be satisfied:

(1) $E$ \ is an UMD space, $p\in \left( 1,\infty \right) $ and $0<h\leq
h_{0}<\infty $ are certain parameters;

(2) \ $m$ is a positive integer and $\alpha =\left( \alpha _{1},\alpha
_{2},...,\alpha _{n}\right) $ are $n$-tuples of nonnegative integer such
that 
\begin{equation*}
\varkappa =\frac{\left\vert \alpha \right\vert }{m}\leq 1,\text{ }0\leq \mu
\leq 1-\varkappa ;
\end{equation*}

(3) \ $A$ is an $R$-sectorial operator in $E$ with\ $0\leq \varphi <\pi $;

Then, the embedding $D^{\alpha }W^{m,p,\gamma }\left( \mathbb{R}%
_{+}^{n};E\left( A\right) ,E\right) \subset L^{p,\lambda }\left( \mathbb{R}%
_{+}^{n};E\left( A^{1-\varkappa -\mu }\right) \right) $ is continuous and
there exists a positive constant $C_{\mu }$ such that the following uniform
estimate holds for all $h\in \left( 0,h_{0}\left. {}\right] \right. $

\begin{equation*}
\left\Vert D^{\alpha }u\right\Vert _{L^{p,\lambda }\left( \mathbb{R}%
_{+}^{n};E\left( A^{1-\varkappa -\mu }\right) \right) }\leq
\end{equation*}

\begin{equation*}
C_{\mu }\left[ h^{\mu }\left\Vert u\right\Vert _{W^{m,p,\lambda }\left( 
\mathbb{R}_{+}^{n};E\left( A\right) ,E\right) }+h^{-\left( 1-\mu \right)
}\left\Vert u\right\Vert _{L^{p,\gamma }\left( \mathbb{R}_{+}^{n};E\right) }%
\right] .
\end{equation*}

Consider the the problem $\left( 1.1\right) -\left( 1.2\right) .$

The function $u\in W^{2,p}\left( \mathbb{R}_{+}^{n};E\left( A\right)
,E\right) $ satisfying the equation $\left( 1.1\right) $ a.e. on $\mathbb{R}%
_{+}^{n}$ is said to be solution of the problem $\left( 1.1\right) .$

\bigskip\ Let $\omega _{k1}=\omega _{k1}\left( x\right) $, $\omega
_{k2}=\omega _{k2}\left( x\right) $ be roots of equations 
\begin{equation}
a_{k}\left( x\right) \omega ^{2}+1=0\text{, }k=1,2,...,n.  \tag{1.10}
\end{equation}
Let $\eta _{1},\eta _{2},...,\eta _{n}$, $\eta _{0}$ denotes $VMO$ modulus
of $a_{1},a_{2,}...,a_{n}$ and $A\left( .\right) A^{-1}\left( x_{0}\right) $%
, respectively.

Let $Q_{0}$ denote the operator in $L^{p}\left( \mathbb{R}_{+}^{n};E\right) $
generated by problem $\left( 3.1\right) $ with $A_{k}=0$ and $\mu =0$, i.e. 
\begin{equation*}
D\left( Q_{0}\right) =W^{2,p}\left( \mathbb{R}_{+}^{n};E\left( A\right)
,E,L_{1}\right) ,\text{ }Q_{0}u=\Sigma_{k=1}^{n}a_{k}\left( x_{k}\right) 
\frac{\partial ^{2}u}{\partial x_{k}^{2}}+A\left( x\right) u.
\end{equation*}

\textbf{Condition 1. }Assume the following conditions are satisfied:

(1) $E$ is an $UMD$ space, $A\left( x\right) $ is an uniformly $R$-positive
operator in $E$ and $A\left( .\right) A^{-1}\left( x_{0}\right) \in
L_{\infty }\left( \mathbb{R}_{+}^{n};L\left( E\right) \right) \cap VMO\left(
L\left( E\right) \right) $ for some $x_{0}\in \mathbb{R}_{+}^{n};$

(2) $a_{k}\in VMO\cap L^{\infty }\left( \mathbb{R}\right) ,$ $-a_{k}\left(
x_{k}\right) \in S\left( \varphi \right) ,$ $a_{k}\left( x_{k}\right) $ $%
\neq 0$ for a.e. $x_{k}\in \mathbb{R}$ and $\alpha _{i}\neq 0$, $i=0,1;$

(3) $\func{Re}\omega _{ki}\neq 0$ and $\frac{\lambda }{\omega _{ki}}\in
S\left( \varphi \right) $ for $\lambda \in S\left( \varphi \right) $, for
a.e. $x_{k}\in \mathbb{R},$ $0\leq \varphi <\pi ,$ $i=1,2$, $k=1,2,...,n$
and $p\in \left( 1,\infty \right) $.

From $\left[ \text{27, Theorems 6.1}\right] $ we have the following result:

\textbf{Theorem A}$_{2}$\textbf{. }Assume the Condition 1.1 holds. Then:

(a) for all $f\in L^{p}\left( \mathbb{R}_{+}^{n};E\right) $, $\mu \in
S\left( \varphi \right) $ and for large enough $\left\vert \lambda
\right\vert ,$ problem $\left( 3.1\right) $ has a unique solution $u\in
W^{2,p}\left( \mathbb{R}_{+}^{n};E\left( A\right) ,E\right) .$ Moreover, the
following coercive uniform estimate holds%
\begin{equation*}
\sum\limits_{k=1}^{n}\sum\limits_{i=0}^{2}\left\vert \mu \right\vert ^{1-%
\frac{i}{2}}\left\Vert \frac{\partial ^{i}u}{\partial x_{k}^{i}}\right\Vert
_{L^{p}\left( \mathbb{R}_{+}^{n};E\right) }+\left\Vert Au\right\Vert
_{L^{p}\left( \mathbb{R}_{+}^{n};E\right) }\leq C\left\Vert \ f\right\Vert
_{L^{p}\left( \mathbb{R}_{+}^{n};E\right) };
\end{equation*}

(b) the operator $Q_{0}$ is $R-$positive in $L^{p}\left( \mathbb{R}%
_{+}^{n};E\right) .$

\bigskip Let 
\begin{equation*}
\Omega =\mathbb{R}_{+}^{n}\times \left( 0,T\right) \text{, }L^{\mathbf{p}%
}\left( \mathbb{\Omega };E\right) =Y\text{, }
\end{equation*}

\begin{equation*}
Y^{2,1,\mathbf{p}}\left( A\right) =W^{2,1,\mathbf{p}}\left( \mathbb{\Omega }%
;E\left( A\right) ,E\right) \text{, }\mathbf{p=}\left( p,p_{1}\right) ,
\end{equation*}%
where 
\begin{equation*}
Y^{2,1,\mathbf{p}}\left( A\right) =\left\{ {}\right. u\in Y,\text{ }\frac{%
\partial u}{\partial t}\in Y,\text{ }\frac{\partial ^{2}u}{\partial x_{k}^{2}%
}\in Y,
\end{equation*}%
\begin{equation*}
\left\Vert u\right\Vert _{Y^{2,1,\mathbf{p}}\left( A\right) }+\left\Vert 
\frac{\partial u}{\partial t}\right\Vert
_{Y}+\sum\limits_{k=1}^{n}\left\Vert \frac{\partial ^{2}u}{\partial x_{k}^{2}%
}\right\Vert _{Y}+\left\Vert Au\right\Vert _{Y}<\infty \left. {}\right\} .
\end{equation*}

From $\left[ \text{27, Theorems 7.1}\right] $ we have

\textbf{Theorem A}$_{3}$\textbf{. }Assume the Condition 1.1 holds for $%
\varphi >\frac{\pi }{2}$. Then for all $f\in L^{p}\left( \mathbb{\Omega }%
;E\right) $ problem $\left( 1.4\right) $ has a unique solution $u\in Y^{2,1,%
\mathbf{p}}\left( A\right) .$ Moreover, the following coercive estimate holds%
\begin{equation*}
\left\Vert \frac{\partial u}{\partial t}\right\Vert
_{Y}+\sum\limits_{k=1}^{n}\left\Vert \frac{\partial ^{2}u}{\partial x_{k}^{2}%
}\right\Vert _{Y}+\left\Vert Au\right\Vert _{Y}\leq C\left\Vert \
f\right\Vert _{Y}.
\end{equation*}

\bigskip We need the abstract versions of Sobolev and Sobolev--Poincare%
%TCIMACRO{\U{b4} }%
%BeginExpansion
\'{}
%EndExpansion
inequalities. Consider the abstract Sobolev space $W_{0}^{m,p}\left( \Omega
;E\right) $ defined the following seminorm 
\begin{equation*}
\left\Vert u\right\Vert _{W_{0}^{m,p}\left( \Omega ;E\right)
}=\sum\limits_{\left\vert \alpha \right\vert =m}\left\Vert D^{\alpha
}u\right\Vert _{L^{p}\left( \Omega ;E\right) }<\infty .
\end{equation*}

Let $f\in L_{loc}(\mathbb{R}^{n};E)$. For a subset $S\subset \mathbb{R}^{n}$
with nonzero Lebesgue measure denoted by $\left\vert S\right\vert $. We
define \ 
\begin{equation*}
\oint_{S}f\left( x\right) dx=\frac{1}{\left\vert S\right\vert }%
\int\limits_{S}f\left( x\right) dx\text{, \ }f_{B(x,r)}=\frac{1}{\left\vert
B(x,r)\right\vert }\int\limits_{B(x,r)}f\left( y\right) dy.
\end{equation*}%
By resoning as in $\left[ \text{31, \S\ 2.4}\right] $ and $\left[ \text{12, 
\S\ 4.5}\right] $ we have the following results

\textbf{Theorem B}$_{1}$\textbf{( Sobolev inequality)}. Assume $E$ is an UMD
space.\ Let $n>mp$ and $q=np\left( n-mp\right) ^{-1}$. Then the following
embedding holds 
\begin{equation*}
W_{0}^{m,p}\left( \mathbb{R}^{n};E\right) \subset L^{q}\left( \mathbb{R}%
^{n};E\right) ,
\end{equation*}

\bigskip i.e. there exists a finite constant $C$ such that for every $u\in
C_{0}^{\infty }\left( \mathbb{R}^{n};E\right) $ the following estimate holds 
\begin{equation*}
\left\Vert u\right\Vert _{L^{q}\left( \mathbb{R}^{n};E\right) }\leq
C\left\Vert u\right\Vert _{W_{0}^{m,p}\left( \mathbb{R}^{n};E\right) }.
\end{equation*}

\textbf{Theorem B}$_{2}$\textbf{( Sobolev--Poincare%
%TCIMACRO{\U{b4} }%
%BeginExpansion
\'{}
%EndExpansion
inequality)}. Assume $E$ is an UMD space. \ Let $n>mp$ and $q=np\left(
n-mp\right) ^{-1}$. Then there exists a finite constant $K$ such that for
every $u\in C_{0}^{\infty }\left( \mathbb{R}^{n};E\right) $ the following
estimate holds 
\begin{equation*}
\left\Vert u-u_{B\left( x,r\right) }\right\Vert _{L^{q}\left( \mathbb{R}%
^{n};E\right) }\leq K\left\Vert u\right\Vert _{W_{0}^{m,p}\left( \mathbb{R}%
^{n};E\right) }.
\end{equation*}

\begin{center}
\bigskip \textbf{2. Abstract elliptic equations in half spaces}
\end{center}

In this section we consider the problem $\left( 1.1\right) -\left(
1.2\right) $ in Morrey space $L^{p,\lambda }\left( \mathbb{R}%
_{+}^{n};E\right) .$

The unction $u\in W^{2,p,\lambda }\left( \mathbb{R}_{+}^{n};E\left( A\right)
,E\right) $ satisfying the equation $\left( 1.1\right) $ a.e. on $\mathbb{R}%
_{+}^{n}$ is said to be solution of the problem $\left( 1.1\right) $ in
Morrey space $L^{p,\lambda }\left( \mathbb{R}_{+}^{n};E\right) .$

\bigskip Let $\eta _{1},\eta _{2},...,\eta _{n}$, $\eta _{0}$ denotes $VMO$
modulus of $a_{1},a_{2,}...,a_{n}$ and $A\left( .\right) A^{-1}\left(
x_{0}\right) $, respectively.

Let $Q_{0}$ denote the operator in $L^{p,\lambda }\left( \mathbb{R}%
_{+}^{n};E\right) $ generated by problem $\left( 3.1\right) $ with $A_{k}=0$
and $\mu =0$, i.e. 
\begin{equation}
D\left( Q_{0}\right) =W^{2,p,\lambda }\left( \mathbb{R}_{+}^{n};E\left(
A\right) ,E,L_{1}\right) ,\text{ }Q_{0}u=\Sigma_{k=1}^{n}a_{k}\left(
x_{k}\right) \frac{\partial ^{2}u}{\partial x_{k}^{2}}+A\left( x\right) u. 
\tag{2.1}
\end{equation}

In this section we will prove the following result:

\textbf{Theorem 2.1. }Assume the Condition 1 holds. Then for all $f\in
L^{p,\lambda }\left( \mathbb{R}_{+}^{n};E\right) $, $\mu \in S\left( \varphi
\right) $ and for large enough $\left\vert \lambda \right\vert ,$ problem $%
\left( 1.1\right) $ has a unique solution $u\in W^{2,p,\lambda }\left( 
\mathbb{R}_{+}^{n};E\left( A\right) ,E\right) .$ Moreover, the following
coercive uniform estimate holds%
\begin{equation}
\sum\limits_{k=1}^{n}\sum\limits_{i=0}^{2}\left\vert \mu \right\vert ^{1-%
\frac{i}{2}}\left\Vert \frac{\partial ^{i}u}{\partial x_{k}^{i}}\right\Vert
_{L^{p,\lambda }\left( \mathbb{R}_{+}^{n};E\right) }+\left\Vert
Au\right\Vert _{L^{p,\lambda }\left( \mathbb{R}_{+}^{n};E\right) }\leq
C\left\Vert \ f\right\Vert _{L^{p,\lambda }\left( \mathbb{R}%
_{+}^{n};E\right) };  \tag{2.2}
\end{equation}

To derive the coercive estimates $\left( 2.2\right) $ in $E$-valued abstact
Morrey spaces, we use the $L^{p}$ estimates in a particular form. This
section is devoted to stating and deriving these estimates from known
results for operators with VMO coefficients. For proving this theorem we
need some Sobolev type inequalities in abstract Sobolev spaces.

For simplicity we will denote the ball of center $0$ and radius $r$ by $%
B\left( r\right) $. Consider a Lipschitz function $\Phi =\left( \Phi
_{1},\Phi _{2},...,\Phi _{n}\right) $: $B\left( 1\right) \rightarrow \mathbb{%
R}^{n}$, and we set $\Phi _{B}=\Phi \left( B\left( 1\right) \right) $, such
that this function has a inverse $\Phi ^{-1}$ with Lipschitz inverse. We
write $J\left( \Phi \right) $ for the Jacobian determinant%
\begin{equation*}
J\left( \Phi \right) =\det \left[ \frac{\partial \Phi _{k}}{\partial x_{j}}%
\right] \text{, }k\text{, }j=1,2,...,n;
\end{equation*}%
and we set 
\begin{equation*}
\Phi _{0}=\left\Vert J\left( \Phi \right) \right\Vert _{L^{\infty }\left(
B\left( 1\right) \right) }+\left\Vert J\left( \Phi ^{-1}\right) \right\Vert
_{L^{\infty }\left( B\left( 1\right) \right) }.
\end{equation*}

Finally, we define 
\begin{equation*}
\Phi _{B}\left( r\right) =\left\{ x:\frac{x}{r}\in \Phi _{B}\right\} \text{
for }r>0.
\end{equation*}

By reasoning as in $\left[ \text{14, Lemma 1}\right] $ and using Theorem B$%
_{2}$, we obtain the following simple variant of the usual Sobolev--Poincare%
%TCIMACRO{\U{b4} }%
%BeginExpansion
\'{}
%EndExpansion
inequality.

\textbf{Lemma 2.1. }Let $E$ be an $UMD$ space, $1\leq p<n$, $q=q=np\left(
n-mp\right) ^{-1}$ and $r>0$. Moreover, assume $u\in W^{1,p}\left( \Phi
_{B}\left( r\right) ;E\right) $ and suppose\ also that 
\begin{equation*}
\int\limits_{\Phi _{B}\left( r\right) }u\left( x\right) dx=0.
\end{equation*}

\bigskip Then there is a constant $C$, determined only by $n$, $p$, $\Phi
_{0}$ and $E$ such that 
\begin{equation}
\left\Vert u\right\Vert _{L^{q}\left( \Phi _{B}\left( r\right) ;E\right)
}\leq C\left\Vert \nabla u\right\Vert _{L^{p}\left( \Phi _{B}\left( r\right)
;E\right) }.  \tag{2.3}
\end{equation}

\textbf{Lemma 2.2. }Let $E$ be an $UMD$ space, $1\leq p<n$, $q=q=np\left(
n-mp\right) ^{-1}$, $r>0$ and $u\in W^{1,p}\left( \Phi _{B}\left( r\right)
;E\right) .$ Assume aso that there are constants $\sigma ,\tau \in \left(
0,1\right) $ and a point $x_{0}\in \partial \Phi _{B}\left( r\right) $ such
that $\left\vert B\left( x_{0},\tau r\right) \cap \Phi _{B}\left( r\right)
\right\vert \geq \sigma r^{n}$, $u=0$ on $B\left( x_{0},\tau r\right) \cap
\Phi _{B}\left( r\right) $, and every line parallel to the $x_{n}$-axis
which intersects $B\left( x_{0},\tau r\right) \cap \Phi _{B}\left( r\right) $
does so at one point. Then there is a constant $C$, determined only by $n$, $%
p$, $\Phi _{0}$, $E$ and $\tau $ such that $\left( 2.3\right) $ holds.

\textbf{Proof. }Let $\Delta =B\left( x_{0},\tau r\right) \cap \Phi
_{B}\left( r\right) $, $w\in L\left( \Delta ;E\right) $ and $q=\frac{np}{n-p}%
.$ Here, 
\begin{equation*}
w\left( \Delta \right) = \oint_{\Delta }w\left( x\right) dx.
\end{equation*}

By reasoning as in $\left[ \text{15, Lemma 2.2}\right] $, for $u\in
L^{p}\left( \Delta ;E\right) $ \ we show that 
\begin{equation*}
\left\Vert u\left( \Delta \right) \right\Vert _{E}\leq C\left( n,p,\sigma
,E\right) r^{\frac{p-n}{p}}\left\Vert \nabla u\right\Vert _{L^{p}\left(
\Delta ;E\right) }.
\end{equation*}

So H\"{o}lder's inequality yields 
\begin{equation}
\left\Vert u\left( \Delta \right) \right\Vert _{L^{q}\left( \Phi _{B}\left(
r\right) ;E\right) }\leq C\left( n,p,\sigma ,E\right) \left\Vert \nabla
u\right\Vert _{L^{p}\left( \Phi _{B}\left( r\right) ;E\right) }.  \tag{2.4}
\end{equation}

\bigskip We observe that 
\begin{equation*}
u\left( \Delta \right) =\frac{1}{\left\vert \Phi _{B}\left( r\right)
\right\vert }\int\limits_{\Phi _{B}\left( r\right) }u\left( x\right) f\left(
x\right) dx
\end{equation*}

with $f\left( x\right) =\frac{\left\vert \Phi _{B}\left( r\right)
\right\vert }{\left\vert \Delta \right\vert }\chi _{\Delta }$, where $\chi
_{\Delta }$ denotes the characteristic function of $\Delta $. Then a similar
way as in $\left[ \text{16, Lemma 6.13}\right] $ with $f$ and $\eta \equiv 
\frac{1}{\left\vert B\left( 1\right) \right\vert }$\ we get that 
\begin{equation}
\left\Vert u-u\left( \Delta \right) \right\Vert _{L^{q}\left( \Phi
_{B}\left( r\right) ;E\right) }\leq C\left( n,p,\sigma ,E\right) \left\Vert
u-u\left( \Phi _{B}\left( r\right) \right) \right\Vert _{L^{p}\left( \Phi
_{B}\left( r\right) ;E\right) }.  \tag{2.5}
\end{equation}

So, Lemma 2.1 applied to $\upsilon =u-u\left( \Delta \right) $ implies that 
\begin{equation}
\left\Vert u-u\left( \Phi _{B}\left( r\right) \right) \right\Vert
_{L^{q}\left( \Phi _{B}\left( r\right) ;E\right) }\leq C\left( n,p,\sigma
,E\right) \left\Vert \nabla u\right\Vert _{L^{p}\left( \Phi _{B}\left(
r\right) ;E\right) }.  \tag{2.6}
\end{equation}

If $p<\frac{n}{2},$ then we also have

Then the inequalities $\left( 2.4\right) -\left( 2.6\right) $ and
Minkowski's inequality gives us $\left( 2.3\right) .$

\textbf{Remark 2.1. }If $u$ vanishes on all of $\partial \Phi _{B}\left(
r\right) $, then Theorem B$_{1}$ implies $(2.4)$ with $C$ determined only by 
$n$, $p$ and $E$.

\bigskip \textbf{Lemma 2.3. }Let $E$ be an $UMD$ space, $1\leq p<n$, $%
q=np\left( n-p\right) ^{-1}$, $r>0$ and $u\in W^{2,p}\left( \Phi _{B}\left(
r\right) ;E\right) .$ If 
\begin{equation*}
\int\limits_{\Phi _{B}\left( r\right) }u\left( x\right) dx=0\text{, }%
\int\limits_{\Phi _{B}\left( r\right) }\nabla u\left( x\right) dx,
\end{equation*}

then there is a constant $C$, determined only by $n$, $p$, $\Phi _{0}$ and $%
E $ such that%
\begin{equation}
\frac{1}{r}\left\Vert u\right\Vert _{L^{q}\left( \Phi _{B}\left( r\right)
;E\right) }+\left\Vert \nabla u\right\Vert _{L^{q}\left( \Phi _{B}\left(
r\right) ;E\right) }\leq C\left( n,p,\sigma ,E\right) \left\Vert \nabla
^{2}u\right\Vert _{L^{p}\left( \Phi _{B}\left( r\right) ;E\right) }. 
\tag{2.7}
\end{equation}

\.{I}f $p<\frac{n}{2}$, then we also have

\begin{equation}
\left\Vert u\right\Vert _{L^{q_{1}}\left( \Phi _{B}\left( r\right) ;E\right)
}\leq C\left( n,p,\sigma ,E\right) \left\Vert \nabla ^{2}u\right\Vert
_{L^{p}\left( \Phi _{B}\left( r\right) ;E\right) },  \tag{2.8}
\end{equation}

with%
\begin{equation*}
q_{1}=\frac{np}{n-2p}.
\end{equation*}

\textbf{Proof. }The estimate for $\left\Vert \nabla u\right\Vert
_{L^{q}\left( \Phi _{B}\left( r\right) ;E\right) }$ follows from Lemma 2.1
and the estimate on $\left\Vert u\right\Vert _{L^{q}\left( \Phi _{B}\left(
r\right) ;E\right) }$ in $(2.7)$ follows from Theorem B$_{2}.$ If $p<\frac{n%
}{2}$, then $q<n$ and Lemma 3.1 yields $\left( 2.8\right) .$

\textbf{Lemma 2.4. }Let $E$ be an $UMD$ space, $1\leq p<n$, $q=np\left(
n-mp\right) ^{-1}$, $r>0$ and $u\in W^{2,p}\left( \Phi _{B}\left( r\right)
;E\right) .$ Suppose also that there are constants $\sigma ,\tau \in \left(
0,1\right) $ and a point $x_{0}\in \partial \Phi _{B}\left( r\right) $ such
that $\left\vert B\left( x_{0},\tau r\right) \cap \Phi _{B}\left( r\right)
\right\vert \geq \sigma r^{n}$, $u=0$ and $\nabla u=0$ on $B\left(
x_{0},\tau r\right) \cap \Phi _{B}\left( r\right) $, and every line parallel
to the $x_{n}$-axis which intersects $B\left( x_{0},\tau r\right) \cap \Phi
_{B}\left( r\right) $ does so at one point. Then there is a constant $C$,
determined only by $n$, $p$, $\Phi _{0}$, $E$ and $\tau $ such that $\left(
2.3\right) $ holds. If $p<\frac{n}{2}$, then $\left( 2.8\right) $ holds with 
$C$\ determined only by $n$, $p$, $\Phi _{0}$, $\sigma $ and $E.$

\textbf{Proof. }We argue as in Lemma 2.3 but using Lemma 2.2 in place of
Lemma 2.1.

\bigskip For our next study, let 
\begin{equation*}
B^{+}\left( r\right) =\left\{ x\text{: }x\in B\left( r\right) \text{, }%
x_{n}>0\right\} ,\text{ }\Phi _{B}^{0}\left( r\right) =\left\{ x\text{: }%
x\in \partial \Phi _{B}\left( r\right) \text{, }x_{n}=0\right\} \text{, }
\end{equation*}%
\begin{equation*}
B^{0}\left( r\right) =\left\{ x\text{: }x\in B\left( r\right) \text{, }%
x_{n}=0\right\} .
\end{equation*}

\bigskip \textbf{Corollary 2.1. }Let $E$ be an $UMD$ space, $1\leq p<n$, $%
q=np\left( n-mp\right) ^{-1}$, $r>0$ and assume%
\begin{equation*}
B^{+}\left( \frac{1}{2}\right) \subset \Phi _{B}\left( 1\right) \subset
B^{+}\left( 1\right) .
\end{equation*}%
If $u\in W^{2,p}\left( \Phi _{B}\left( r\right) ;E\right) ,$ $u=0$ on $\Phi
_{B}^{0}\left( r\right) ,$ and $\int\limits_{\Phi _{B}\left( r\right) }\frac{%
\partial }{\partial x_{n}}u\left( x\right) dx=0$, then $\left( 2.7\right) $
holds. If $p<\frac{n}{2}$, then $\left( 2.8\right) $ holds with $C$\
determined only by $n$, $p$, $\Phi _{0}$, $\sigma $ and $E.$

\textbf{Proof. }We apply \ here, Lemma 2.2 to $\frac{\partial u}{\partial
x_{1}}$, $\frac{\partial u}{\partial x_{2}}$, ...,$\frac{\partial u}{%
\partial x_{n-1}}$ and Lemma 2.1 to $\frac{\partial u}{\partial x_{n}}$ to
deduce the estimate on $\nabla u$. From this estimate, we infer the estimate
on $u$ by applying Lemma 2.2 and by applying Theorem B$_{2}$ in general.

\textbf{Corollary 2.2. }Let $E$ be an $UMD$ space, $1\leq p<n$, $q=np\left(
n-mp\right) ^{-1}$, $r>0$. Assume 
\begin{equation*}
B^{+}\left( \frac{1}{2}\right) \subset \Phi _{B}\left( 1\right) \subset
B^{+}\left( 1\right) .
\end{equation*}%
Suppose $\beta =\left( \beta _{1}\text{, }\beta _{2}\text{, ..., }\beta
_{n}\right) ,$ where $\beta _{i}$ are constant with $\beta _{n}>0$ that $%
\beta .\nabla u=0$ on $\Phi _{B}^{0}\left( r\right) $. Moreover, for $u\in
W^{2,p}\left( \Phi _{B}\left( r\right) ;E\right) ,$ 
\begin{equation*}
\int\limits_{\Phi _{B}\left( r\right) }\frac{\partial }{\partial x_{n}}%
u\left( x\right) dx=0\text{, }\int\limits_{\Phi _{B}\left( r\right) }\left[
\nabla u-\left( \beta .\nabla u\right) \beta \right] dx=0,
\end{equation*}
then $\left( 2.7\right) $ holds, if $p<\frac{n}{2}$, then $\left( 2.8\right) 
$ holds with $C$\ determined also by $\chi .$

\textbf{Proof. }This time we use Lemma 2.1 and on the components of $\nabla
u $ which are perpendicular to $\beta $ and Lemma 2.2 on $\beta .\nabla u.$

Consider the problem 
\begin{equation}
\Sigma_{k=1}^{n}a_{k}\left( x_{k}\right) \frac{\partial ^{2}u}{\partial
x_{k}^{2}}+A\left( x\right) u=f\left( x\right) ,\text{ }x\in B\left(
r\right) ,  \tag{2.9}
\end{equation}

\begin{equation*}
L_{1}u=\sum\limits_{i=0}^{\sigma }\alpha _{i}u^{\left( i\right) }\left(
x\right) =0\text{, for }x\in \partial B,
\end{equation*}%
where $\sigma \in \left\{ 0,1\right\} $, $a_{k}$ are complex-valued
functions and $A=A\left( x\right) $ is a linear operator in a Banach space $%
E.$ Let 
\begin{equation*}
Q_{0}u=\Sigma _{k=1}^{n}a_{k}\left( x_{k}\right) \frac{\partial ^{2}u}{%
\partial x_{k}^{2}}+A\left( x\right) u\text{ for }W^{2,p}\left( B\left(
r\right) ;E\left( A\right) ,E\right) .
\end{equation*}

Let $\eta _{1},\eta _{2},...,\eta _{n}$, $\eta _{0}$ denotes $VMO$ modulus
of $a_{1},a_{2,}...,a_{n}$, $A\left( .\right) A^{-1}\left( x_{0}\right) $,
respectively and $\eta =\left( \eta _{1},\eta _{2},...,\eta _{n},\eta
_{0}\right) .$ By resoning as in the proof of Theorem A$_{2}$ we obtain the
following result

\textbf{Theorem 2.2. }Assume the Condition 1 holds. Let $p\in \left(
1,\infty \right) $ and $f\in L^{p}\left( B\left( r\right) ;E\right) $. Then
there is a unique solution%
\begin{equation*}
\upsilon \in W_{0}^{1,p}\left( B\left( r\right) ;E\right) \cap W^{2,p}\left(
B\left( r\right) ;E\left( A\right) ,E\right)
\end{equation*}
of the problem $\left( 2.9\right) $\ such that 
\begin{equation}
\left\Vert \upsilon \right\Vert _{W^{2,p}\left( B\left( r\right) ;E\left(
A\right) ,E\right) }\leq C\left( n,\eta ,p,\lambda ,E,A\right) \left\Vert \
f\right\Vert _{L^{p}\left( B\left( r\right) ;,E\right) }.  \tag{2.10}
\end{equation}

\bigskip \textbf{Remark 2.2. }By $\left[ \text{23, Theorem 4.2}\right] $ we
get that for all $u\in W^{2,p}\left( B\left( r\right) ;E\right) $ or $u\in $ 
$W^{2,p}\left( B\left( r\right) ;E\left( A\right) ,E\right) $ we have,
respectively,%
\begin{equation*}
\left\Vert u\right\Vert _{W^{2,p}\left( B\left( r\right) ;E\right) }\sim
\left\Vert \nabla ^{2}u\right\Vert _{L^{p}\left( B\left( r\right) ;E\right)
}+\left\Vert u\right\Vert _{L^{p}\left( B\left( r\right) ;E\right) }
\end{equation*}
\begin{equation*}
\left\Vert u\right\Vert _{W^{2,p}\left( B\left( r\right) ;E\left( A\right)
,E\right) }\sim \left\Vert \nabla ^{2}u\right\Vert _{L^{p}\left( B\left(
r\right) ;E\right) }+\left\Vert Au\right\Vert _{L^{p}\left( B\left( r\right)
;E\right) }.
\end{equation*}

By following $\left[ \text{16, Theorem 8}\right] $, by using Lemmas 2.1,
2.3, 2.4, Remark 2.2 and Theorem 2.2 we obtain the following

\bigskip \textbf{Theorem 2.3. }Assume the Condition 1 holds. Let $1\leq
p<q<\infty $ and suppose that $w\in W^{2,p}\left( B\left( r\right)
;,E\right) $ satisfies $Q_{0}w=0$ in $B\left( r\right) .$Then $w\in
W^{2,p}\left( B\left( \frac{r}{2}\right) ;E\right) $ \ and, if $r\leq 1$,
then there is $C\left( \eta ,p,q,\lambda ,E\right) $ such that 
\begin{equation}
\left\Vert \nabla ^{2}w\right\Vert _{L^{q}\left( B\left( \frac{r}{2}\right)
;E\right) }\leq C\left( \eta ,p,q,\lambda ,E,A\right) r^{n\left( \frac{1}{q}-%
\frac{1}{p}\right) }\left\Vert \ \nabla ^{2}w\right\Vert _{L^{p}\left(
B\left( r\right) ;E\right) }.  \tag{2.11}
\end{equation}

\bigskip \textbf{Proof. }\.{I}ndeed, let $\zeta $ be a nonnegative,
compactly supported $C^{\left( 2\right) }\left( B\left( r\right) \right) $
function with $\zeta \equiv 1$ in $B\left( \frac{r}{2}\right) $ and $%
r\left\vert \nabla \zeta \right\vert +r^{2}\left\vert \nabla ^{2}\zeta
\right\vert \leq C\left( n\right) $. We then define

\begin{equation*}
\bar{w}=w\left( x\right) -x\oint_{B\left( r\right) }\nabla w\left( x\right)
dx- \oint_{B\left( r\right) }w\left( x\right) dx
\end{equation*}%
and set $w_{k}=\zeta ^{k}\bar{w}$. First we show 
\begin{equation*}
\left\Vert \nabla ^{2}w_{k}\right\Vert _{L^{q\left( k\right) }\left(
E\right) }\leq Cr^{-k}\left\Vert \nabla ^{2}w\right\Vert _{L^{1}\left(
E\right) },
\end{equation*}%
here 
\begin{equation*}
q\left( k\right) =n\left( n-k\right) ^{-1},k=0,1,2,...n-1.
\end{equation*}

Then by reasoning as in $\left[ \text{16, Theorem 8}\right] $ we obtain the
assertion.

\bigskip We can now, give the interier Morrey space esimates for solution of
abstract elliptic equation $\left( 2.9\right) .$ By following $\left[ \text{%
16, Theorem 9}\right] $ we show the following:

\textbf{Theorem 2.4. }Let the assumptions (1)-(3) of Condition 1 be
satisfied. Assume $u\in W^{2,p}\left( B\left( r\right) ;E\left( A\right)
,E\right) $ is a solution of the problem $\left( 2.9\right) $ in $%
L^{p,\lambda }\left( B\left( r\right) ;E\right) $. Then $\nabla ^{2}u\in
L^{p,\lambda }\left( B\left( \frac{r}{2}\right) ;E\right) $ and there is a
constant $C=$ $C\left( \eta ,p,r,\lambda ,E,A\right) $ such that 
\begin{equation}
\left\Vert \nabla ^{2}u\right\Vert _{L^{p,\lambda }\left( B\left( \frac{r}{2}%
\right) ;E\right) }+\left\Vert Au\right\Vert _{L^{p,\lambda }\left( B\left( 
\frac{r}{2}\right) ;E\right) }\leq  \tag{2.12}
\end{equation}%
\begin{equation*}
C\left( \left\Vert Q_{0}u\right\Vert _{L^{p,\lambda }\left( B\left( r\right)
;E\right) }+\left\Vert \ \nabla ^{2}u\right\Vert _{L^{p}\left( B\left(
r\right) ;E\right) }+\left\Vert Au\right\Vert _{L^{p}\left( B\left( r\right)
;E\right) }\right) .
\end{equation*}

\textbf{Proof. }Let $x_{0}\in B\left( \frac{r}{2}\right) $ and let $\rho \in
\left( 0.r-\left\vert x_{0}\right\vert \right) .$ Assume $\upsilon \in
W_{0}^{1,p}\left( B\left( r\right) ;E\right) \cap W^{2,p}\left( B\left(
r\right) ;E\left( A\right) ,E\right) $ is the solution of $\left( 2.9\right) 
$ given by Theorem 2.2. Choose $q>p\left( 1-\frac{\lambda }{n}\right) ^{-1}$
and set $w=u-\upsilon .$ It follows from Theorem 2.3 that 
\begin{equation*}
\left\Vert \nabla ^{2}w\right\Vert _{L^{q}\left( B\left( x_{0},\frac{\rho }{2%
}\right) ;E\right) }\leq C\rho ^{n\left( \frac{1}{q}-\frac{1}{p}\right)
}\left\Vert \ \nabla ^{2}w\right\Vert _{L^{p}\left( B\left( x_{0},\rho
\right) ;E\right) }.
\end{equation*}

Therefore, if $\tau \in \left( 0,\frac{1}{2}\right) $, then by H\"{o}lder's
inequality and $\left( 2.11\right) $\ we have 
\begin{equation*}
\left\Vert \nabla ^{2}w\right\Vert _{L^{p}\left( B\left( x_{0},\tau \rho
\right) ;E\right) }\leq C\left( \tau \rho \right) ^{n\left( \frac{1}{p}-%
\frac{1}{q}\right) }\left\Vert \ \nabla ^{2}w\right\Vert _{L^{q}\left(
B\left( x_{0},\tau \rho \right) ;E\right) }\leq
\end{equation*}%
\begin{equation*}
C\left( \tau \rho \right) ^{n\left( \frac{1}{p}-\frac{1}{q}\right)
}\left\Vert \ \nabla ^{2}w\right\Vert _{L^{q}\left( B\left( x_{0},\frac{\rho 
}{2}\right) ;E\right) }\leq C\tau ^{n\left( \frac{1}{p}-\frac{1}{q}\right)
}\left\Vert \ \nabla ^{2}w\right\Vert _{L^{p}\left( B\left( x_{0},\rho
\right) ;E\right) }.
\end{equation*}

\bigskip In addition, the estimate $(2.10)$, Theorem A$_{2}$ and Remark 2.1
implies that%
\begin{equation*}
\left\Vert \nabla ^{2}\upsilon \right\Vert _{L^{p}\left( B\left( x_{0},\tau
\rho \right) ;E\right) }+\left\Vert A\upsilon \right\Vert _{L^{p}\left(
B\left( x_{0},\tau \rho \right) ;E\right) }\leq
\end{equation*}%
\begin{equation*}
\left\Vert \ \nabla ^{2}\upsilon \right\Vert _{L^{p}\left( B\left(
x_{0},\rho \right) ;E\right) }+\left\Vert \ A\upsilon \right\Vert
_{L^{p}\left( B\left( x_{0},\rho \right) ;E\right) }\leq \left\Vert
f\right\Vert _{L^{p}\left( B\left( x_{0},\rho \right) ;E\right) }\leq
\end{equation*}%
\begin{equation*}
\rho ^{\frac{\lambda }{p}}\left\Vert f\right\Vert _{L^{p}\left( B\left(
x_{0},r\right) ;E\right) }.
\end{equation*}%
Hence, 
\begin{equation*}
\left\Vert \nabla ^{2}u\right\Vert _{L^{p}\left( B\left( x_{0},\tau \rho
\right) ;E\right) }+\left\Vert Au\right\Vert _{L^{p}\left( B\left(
x_{0},\tau \rho \right) ;E\right) }\leq C\left[ \tau ^{n\left( \frac{1}{p}-%
\frac{1}{q}\right) }\right. \left\Vert \ \nabla ^{2}u\right\Vert
_{L^{p}\left( B\left( x_{0},\rho \right) ;E\right) }+
\end{equation*}%
\begin{equation*}
\rho ^{\frac{\lambda }{p}}\left\Vert f\right\Vert _{L^{p}\left( B\left(
x_{0},r\right) ;E\right) }\left. {}\right] .
\end{equation*}%
The standard iteration scheme then implies that%
\begin{equation*}
\left\Vert \nabla ^{2}u\right\Vert _{L^{p}\left( B\left( x_{0},r\right)
;E\right) }+\left\Vert Au\right\Vert _{L^{p}\left( B\left( x_{0},r\right)
;E\right) }\leq
\end{equation*}%
\begin{equation}
C\left( \rho r^{-1}\right) ^{\frac{\lambda }{p}}\left\Vert \ \nabla
^{2}u\right\Vert _{L^{p}\left( B\left( x_{0},\frac{r}{2}\right) ;E\right)
}++\rho ^{\frac{\lambda }{p}}\left\Vert f\right\Vert _{L^{p}\left( B\left(
x_{0},r\right) ;E\right) }.  \tag{2.13}
\end{equation}%
From $\left( 2.13\right) $ we obtain $\left( 2.12\right) .$

Let 
\begin{equation*}
O^{+}\left( r\right) =B\left( r,x\right) \cap \mathbb{R}_{+}^{n}.
\end{equation*}%
Define $W_{L}^{2,p}\left( O^{+}\left( r\right) ;E\left( A\right) ,E\right) $
to be the closure in $W^{2,p}\left( O^{+}\left( r\right) ;E\left( A\right)
,E\right) $ of the subspace $C_{L}\left( E\right) =\left\{ {}\right. u$: $u$
is the restriction to $\mathbb{R}_{+}^{n}$\ of a function belonging to $%
C_{0}^{\infty }\left( B_{r};E\right) $ with $Lu=\sum\limits_{i=0}^{\sigma
}\alpha _{i}u^{\left( i\right) }\left( x^{\prime },0\right) =0$, $x^{\prime
}=\left( x_{1},x_{2},...,x_{n-1}\right) ,\left. \sigma \in \left\{
0,1\right\} \right\} .$

Here,\ we obtain the following variation of Theorem 2.4, i.e. we give
boundary estimates in\ $E$-valued Morrey space for solution of $\left(
2.9\right) $.

\textbf{Theorem 2.5. }Let the Condition 1 holds. Assume $u\in W^{2,p}\left(
B^{+}\left( r\right) ;E\left( A\right) ,E\right) $ and $Q_{0}u\in $ $%
L^{p,\lambda }\left( B\left( r\right) ;E\right) $. If also $L_{1}u=0$ on $%
B^{0}\left( r\right) $, then $\nabla ^{2}u\in L^{p,\lambda }\left( B\left( 
\frac{r}{2}\right) ;E\right) $ and there is a constant $C=$ $C\left( \eta
,p,r,\lambda ,E,A\right) $ such that 
\begin{equation}
\left\Vert \nabla ^{2}u\right\Vert _{L^{p,\lambda }\left( B^{+}\left( \frac{r%
}{2}\right) ;E\right) }+\left\Vert Au\right\Vert _{L^{p,\lambda }\left(
B^{+}\left( \frac{r}{2}\right) ;E\right) }\leq  \tag{2.14}
\end{equation}%
\begin{equation*}
C\left( \left\Vert Q_{0}u\right\Vert _{L^{p,\lambda }\left( B^{+}\left(
r\right) ;E\right) }+\left\Vert \ \nabla ^{2}u\right\Vert _{L^{p}\left(
B^{+}\left( r\right) ;E\right) }+\left\Vert Au\right\Vert _{L^{p}\left(
B^{+}\left( r\right) ;E\right) }\right) .
\end{equation*}

\textbf{Proof. }Let $B^{+}\left( \frac{3}{4}\right) \subset \Phi _{B}\subset
B^{+}\left( 1\right) $ and such that $\Phi _{B}$ is $C^{2}$-diffeomorphic to 
$B\left( 1\right) $. Then for $r>0$ and $x_{0}\in \mathbb{R}^{n}$ with $%
x_{0n}=0$, we define 
\begin{equation*}
\Phi _{B}\left( x_{0},r\right) =\left\{ x\text{: }\left( x-x_{0}\right)
r^{-1}\in \Phi _{B}\right\} .
\end{equation*}

Again we write $\Phi _{B}\left( r\right) $ for $\Phi _{B}\left( 0,r\right) $%
. For the Dirichlet problem we see that Theorem 2.2 holds with $\Phi
_{B}\left( r\right) $ in place of $B\left( r\right) $ with the same proof.
Theorem 2.3 holds with the same changes in the statement but the proof must
be modified somewhat. We take the support of $\zeta $ to be a subset of $%
B\left( \frac{3r}{4}\right) $, and we define $\bar{w}\left( \Delta \right) $
by%
\begin{equation*}
\bar{w}\left( \Delta \right) =w\left( x\right) -x_{n} \oint_{\Delta }\frac{%
\partial }{\partial x_{n}}w\left( x\right) dx.
\end{equation*}%
Then we use the Sobolev inequality Corollary 2.1 in place of Lemma 2.4 and
by reasonins as in $\left[ \text{16, Theorem 10}\right] $, we obtain the
assertion.

\textbf{Theorem 2.6. }Let the Condition 1 holds. Let $r>0$, $p\in \left(
1,\infty \right) $ and $f\in L^{p}\left( \Phi _{B}\left( r\right) ;E\right)
. $\ Then there is $\upsilon \in W^{2,p}\left( \Phi _{B}\left( r\right)
;E\left( A\right) ,E\right) $ such that $Q_{0}\upsilon =f$ \ in $\Phi
_{B}\left( r\right) $ and $Lu=0$ in $\partial \Phi _{B}\left( r\right) .$
Moreover, if $r\leq 1,$ then 
\begin{equation}
\left\Vert \nabla ^{2}\upsilon \right\Vert _{L^{p}\left( \Phi _{B}\left(
r\right) ;E\right) }+\left\Vert A\upsilon \right\Vert _{L^{p}\left( \Phi
_{B}\left( r\right) ;E\right) }\leq  \tag{2.15}
\end{equation}

\begin{equation*}
C\left( n,\eta ,p,\lambda ,E,A\right) \left\Vert f\right\Vert _{L^{p}\left(
\Phi _{B}\left( r\right) ;E\right) }.
\end{equation*}%
\textbf{Proof. }Now the existence and regularity results, as well as the
estimate in $\Phi _{B}\left( r\right) $, is obtained by reasonins as $[$18,
Lemma 3$]$ and $\left[ \text{14, Theorem 3.2}\right] .$

\textbf{Proof of Theorem 2.1. }Indeed, the prove Theorem 2.1 is obtained
from Theorems 2.5, 2.6 by localizasion and perturbation argument.

Consider now the problem $\left( 1.1\right) $. Let $Q$ is operator in $X=$ $%
L^{p,\lambda }\left( \mathbb{R}_{+}^{n};E\right) $ generated by problem $%
\left( 1.1\right) -\left( 1.2\right) $ and $Y=W^{2,p,\lambda }\left( \mathbb{%
R}_{+}^{n};E\left( A\right) ,E\right) .$

Now, we will show the following main theorem of this section

\textbf{Theorem 2.7. }Assume the Condition 1 holds and suppose $A_{j}\left(
.\right) A^{-\left( \frac{1}{2}-\nu \right) }\left( x_{0}\right) $ for $\nu
\in \left( 0,\frac{1}{2}\right) $ and some $x_{0}\in \mathbb{R}_{+}^{n}$.
Then for all $f\in X$, $\mu \in S\left( \varphi \right) $ and for large
enough $\left\vert \lambda \right\vert ,$ problem $\left( 1.1\right) -\left(
1.2\right) $ has a unique solution $u\in Y.$ Moreover, the following
coercive uniform estimate holds%
\begin{equation}
\sum\limits_{k=1}^{n}\sum\limits_{i=0}^{2}\left\vert \mu \right\vert ^{1-%
\frac{i}{2}}\left\Vert \frac{\partial ^{i}u}{\partial x_{k}^{i}}\right\Vert
_{X}+\left\Vert Au\right\Vert _{X}\leq C\left\Vert \ f\right\Vert _{X}. 
\tag{2.16.}
\end{equation}

\textbf{Proof. } It sufficient to show that the operator $Q+\lambda $ has a
bounded inverse\ from $X$ to $Y.$\ Put $Qu=Q_{0}u+Q_{1}u,$ where 
\begin{equation*}
Q_{1}u=\Sigma _{k=1}^{n}A_{k}\left( x\right) \frac{\partial u}{\partial x_{k}%
},\text{ }u\in Y.
\end{equation*}%
By Theorem A$_{1},$ for all $\varepsilon >0$ there is a continuous function $%
C\left( \varepsilon \right) >0$ such that

\begin{equation}
\left\Vert A_{k}\frac{\partial u}{\partial x_{k}}\right\Vert _{X}\leq
C\left\Vert A^{\frac{1}{2}-\mu }\frac{\partial u}{\partial x_{k}}\right\Vert
_{X}\leq \varepsilon \left\Vert u\right\Vert _{Y}+C\left( \varepsilon
\right) \left\Vert u\right\Vert _{X}  \tag{2.17}
\end{equation}%
for $0<\mu <\frac{1}{2}$. \ Since $\left\Vert u\right\Vert _{X}=\frac{1}{%
\left\vert \lambda \right\vert }\left\Vert \left( Q_{0}+\lambda \right)
u-Q_{0}u\right\Vert _{X},$ by definition of the space $X$ and by Theorem
2.1, we have%
\begin{equation*}
\left\Vert u\right\Vert _{X}\leq \frac{1}{\left\vert \lambda \right\vert }%
\left[ \left\Vert \left( Q_{0}+\lambda \right) u\right\Vert _{X}+\left\Vert
Q_{0}u\right\Vert _{X}\right] \leq \frac{C}{\left\vert \lambda \right\vert }%
\left\Vert \left( Q_{0}+\lambda \right) u\right\Vert _{X},
\end{equation*}

In view $\left( 2.17\right) $ then for sufficientli large $\left\vert
\lambda \right\vert $\ again by applying Theorem 2.1 we get%
\begin{equation}
\left\Vert A_{k}\frac{\partial u}{\partial x_{k}}\right\Vert _{L^{\mathbf{p}%
}\left( \mathbb{R}^{n};E\right) }<\delta \left\Vert \left( Q_{0}+\lambda
\right) u\right\Vert _{L^{\mathbf{p}}\left( \mathbb{R}^{n};E\right) }, 
\tag{2.18}
\end{equation}%
\ \ for $u\in Y$ with $\delta <1.$ By Theorem 2.1, the operator $%
Q_{0}+\lambda $ has a bounded inverse $\left( Q_{0}+\lambda \right) ^{-1}$\
from $X$ to $Y$ for sufficiently large $\left\vert \lambda \right\vert .$
So, $\left( 2.2\right) $ and $\left( 2.18\right) $\ implies the following
estimate 
\begin{equation}
\left\Vert Q_{1}\left( Q_{0}+\lambda \right) ^{-1}\right\Vert _{L\left(
X\right) }<1.  \tag{2.19}
\end{equation}%
Then, by $\left( 2.19\right) $ and in view of Theorem 2.1 we get the
assertion.

\begin{center}
\bigskip \textbf{3. Abstract Cauchy problem for parabolic equation in Morrey
spaces}
\end{center}

Consider the mixed problem for abstract parabolic equation $\left(
1.4\right) $ in the Morrey space $L^{p,q,\lambda ,\mu }\left( \Omega
;E\right) ,$ where $\Omega =\mathbb{R}_{+}^{n}\times \left( 0,T\right) $.
Now, let 
\begin{equation*}
Q\left( r\right) =\left\{ Z=\left( x,t\right) \in \mathbb{R}^{n+1}\text{, }%
\left\vert x-x_{0}\right\vert <r\text{, }-r^{2}<t-t_{0}<0\right\} ,
\end{equation*}%
\begin{equation*}
Q^{+}\left( r\right) =\left\{ \left( x,t\right) \in \mathbb{R}^{n+1}\text{, }%
\left\vert x-x_{0}\right\vert <r\text{, }x_{n}>0,-r^{2}<t-t_{0}<0\right\} ,
\end{equation*}%
\begin{equation*}
Q^{0}\left( r\right) =\left\{ \left( x,t\right) \in \mathbb{R}^{n+1}\text{, }%
\left\vert x-x_{0}\right\vert <r\text{, }x_{n}=0,-r^{2}<t-t_{0}<0\right\}
\end{equation*}%
for $Z_{0}=\left( x_{0},t_{0}\right) $\ and $Q\left( r\right) =Q\left(
0,r\right) $. In addition, we define $BMO$ with respect to the abstract
parabolic metric, so we say that $f\in BMO$ if%
\begin{equation*}
\doint\limits_{Q\left( Z_{0},r\right) }\left\Vert f\left( x\right)
-\doint\limits_{Q\left( Z_{0},r\right) }f\left( y\right) dy\right\Vert _{E}dx
\end{equation*}%
is bounded independent of $Z_{0}$ and $r$. \ 

\bigskip Let 
\begin{equation*}
\text{ }L^{p}\left( \mathbb{\Omega };E\right) =Y\text{, }Y^{2,1,p}=W^{2,1,p}%
\left( \mathbb{\Omega };E\right) ,
\end{equation*}%
where 
\begin{equation*}
Y^{2,1,p}=\left\{ {}\right. u\in Y,\text{ }\frac{\partial u}{\partial t}\in
Y,\text{ }\frac{\partial ^{2}u}{\partial x_{k}^{2}}\in Y,\text{ }k=1,2,...,n,
\end{equation*}%
\begin{equation*}
\left\Vert u\right\Vert _{Y^{2,1,p}}+\left\Vert \frac{\partial u}{\partial t}%
\right\Vert _{Y}+\sum\limits_{k=1}^{n}\left\Vert \frac{\partial ^{2}u}{%
\partial x_{k}^{2}}\right\Vert _{Y}<\infty \left. {}\right\} .
\end{equation*}

Here, $\Phi =L^{p,\lambda }\left( \mathbb{\Omega };E\right) $ denotes the
space of all $E$-valued Morrey space. Analogously, $\Phi ^{1,2}\left(
A\right) =W^{1,2,p,\lambda }\left( \Omega ;E\left( A\right) ,E\right) $
denotes the Sobolev-Morrey space with the following mixed norm%
\begin{equation*}
\left\Vert u\right\Vert _{\Phi ^{1,2}\left( A\right) }=\left\Vert \frac{%
\partial u}{\partial t}\right\Vert _{\Phi }+\Sigma _{k=1}^{n}\left\Vert 
\frac{\partial ^{2}u}{\partial x_{k}^{2}}\right\Vert _{\Phi }+\left\Vert
Au\right\Vert _{\Phi }<\infty .
\end{equation*}

We then by resoning as in Lemma 2.3,\ have the following parabolic version
of the Lemma 2.3.

\textbf{Lemma 3.1}. Let $E$ be an $UMD$ space, 
\begin{equation*}
1\leq p<n+2,q=p\left( n+2\right) \left( n+2-p\right) ^{-1},r>0
\end{equation*}
and $u\in W^{2,1,p}\left( Q\left( r\right) ;E\right) .$ If 
\begin{equation*}
\int\limits_{Q\left( r\right) }u\left( x\right) dx=0\text{, }%
\int\limits_{Q\left( r\right) }\nabla u\left( x\right) dx,
\end{equation*}

then there is a constant $C$, determined only by $n$, $p$, $\Phi _{0}$ and $%
E $ such that%
\begin{equation}
\frac{1}{r}\left\Vert u\right\Vert _{L^{q}\left( Q\left( r\right) ;E\right)
}+\left\Vert \nabla u\right\Vert _{L^{q}\left( Q\left( r\right) ;E\right)
}\leq C\left( n,p,\sigma ,E\right) \left\Vert \nabla ^{2}u\right\Vert
_{L^{p}\left( Q\left( r\right) ;E\right) }.  \tag{3.1}
\end{equation}

\.{I}f $p<\frac{n+2}{2}$, then we also have

\begin{equation}
\left\Vert u\right\Vert _{L^{q_{1}}\left( Q\left( r\right) ;E\right) }\leq
C\left( n,p,\sigma ,E\right) \left\Vert \nabla ^{2}u\right\Vert
_{L^{p}\left( Q\left( r\right) ;E\right) },  \tag{3.2}
\end{equation}

with%
\begin{equation*}
q_{1}=p\left( n+2\right) \left( n+2-2p\right) ^{-1}.
\end{equation*}

\bigskip \textbf{Proof.} Indeed, by combine $\left[ \text{9, Lemmas 2, 11}%
\right] $ with $\left[ \text{15, Lemma 2. 3.3}\right] $ and by replacing
Banach space valued function with complex-valued funtion we get the
assertion.

\bigskip Let 
\begin{equation*}
G\left( r\right) =\Phi _{B}\left( r\right) \times \left( -r^{2},0\right) .
\end{equation*}

With the above estimates in hand, we can modify the arguments of Section 1
to prove the parabolic analogs of the \ Lemma 2.4:

\textbf{Lemma 3.2. }Let $E$ be an $UMD$ space and let 
\begin{equation*}
1\leq p<n+2,q=p\left( n+2\right) \left( n+2-p\right) ^{-1},r>0.
\end{equation*}
Assume $u\in W^{2,1,p}\left( G\left( r\right) ;E\right) ,$ $u=0$ and $\nabla
u=0$ on%
\begin{equation*}
B\left( x_{0}r,\tau r\right) \cap \Phi _{B}\left( r\right) \cap \Phi
_{B}\left( r\right) \times \left( -r^{2},0\right) .
\end{equation*}
Suppose also that there are positive constants $\sigma $ and $\tau $ along
with $x_{0}\in \partial \Phi _{B}\left( r\right) $ such that $\left\vert
B\left( x_{0},\tau \right) \cap \partial \Phi _{B}\left( r\right)
\right\vert \geq \sigma $. Moreover, every line parallel to the $x_{n}$-axis
which intersects $B\left( x_{0},\tau \right) \cap \Phi _{B}\left( r\right) $
does so in exactly\ at one point. Then there is a constant $C$, determined
only by $n$, $p$, $\Phi _{0}$, $E$ and $\tau $ such that $\left( 3.1\right) $
holds. If $p<\frac{n+2}{2}$, then $\left( 3.2\right) $ holds with $C$\
determined only by $n$, $p$, $\Phi _{0}$, $\sigma $ and $E.$

For the Cauchy--Dirichlet problem, we have the following analog of Corollary
2.1.

\bigskip \textbf{Corollary 3.1. }Let $E$ be an $UMD$ space and let 
\begin{equation*}
1\leq p<n+2,q=p\left( n+2\right) \left( n+2-p\right) ^{-1},r>0.
\end{equation*}
Assume%
\begin{equation*}
B^{+}\left( \frac{1}{2}\right) \subset \Phi _{B}\left( 1\right) \subset
B^{+}\left( 1\right) .
\end{equation*}%
If $u\in W^{2,1,p}\left( G\left( r\right) ;E\right) ,$ $u=0$ on $\Phi
_{B}^{0}\left( r\right) \times \left( -r^{2},0\right) ,$ and $%
\int\limits_{\Phi _{B}\left( r\right) }\frac{\partial }{\partial x_{n}}%
u\left( x\right) dx=0$, then $\left( 2.7\right) $ holds. If $p<\frac{n+2}{2}$%
, then $\left( 2.8\right) $ holds with $C$\ determined only by $n$, $p$, $%
\Phi _{0}$, $\sigma $ and $E.$

\bigskip In a similar parabolic analog of Corollary 2.2 holds:

\textbf{Corollary 3.2. }Let $E$ be an $UMD$ space and let 
\begin{equation*}
1\leq p<n+2,q=p\left( n+2\right) \left( n+2-p\right) ^{-1},r>0.
\end{equation*}%
Assume 
\begin{equation*}
B^{+}\left( \frac{1}{2}\right) \subset \Phi _{B}\left( 1\right) \subset
B^{+}\left( 1\right) .
\end{equation*}%
Suppose $\beta =\left( \beta _{1}\text{, }\beta _{2}\text{, ..., }\beta
_{n}\right) ,$ where $\beta _{i}$ are constant with $\beta _{n}>0$ that $%
\beta .\nabla u=0$ on $\Phi _{B}^{0}\left( r\right) \times \left(
-r^{2},0\right) $. Moreover, for $u\in W^{2,1,p}\left( G\left( r\right)
;E\right) ,$ 
\begin{equation*}
\int\limits_{\Phi _{B}\left( r\right) }\frac{\partial }{\partial x_{n}}%
u\left( x\right) dx=0\text{, }\int\limits_{\Phi _{B}\left( r\right) }\left[
\nabla u-\left( \beta .\nabla u\right) \beta \right] dx=0,
\end{equation*}%
then $\left( 2.7\right) $ holds, if $p<\frac{n+2}{2}$, then $\left(
2.8\right) $ holds with $C$\ determined also by $\chi .$

\bigskip The analog of bootstrap estimate takes the following form in the
abstract parabolic case.

\textbf{Theorem 3.1. }Assume the Condition 1 holds. Let $1\leq p<q<\infty $
and suppose that $w\in W^{2,p}\left( Q\left( r\right) ;E\right) $ satisfies%
\begin{equation*}
Lw=\frac{\partial w}{\partial t}+\Sigma _{k=1}^{n}a_{k}\left( x_{k}\right) 
\frac{\partial ^{2}w}{\partial x_{k}^{2}}=0\text{ on }Q\left( r\right) .
\end{equation*}
Then $w\in W^{2,p}\left( Q\left( \frac{r}{2}\right) ;E\right) $ \ and, if $%
r\leq 1$, then there is $C\left( \eta ,p,q,\lambda ,E\right) $ such that 
\begin{equation}
\left\Vert w\right\Vert _{W^{2}{}^{,q}\left( Q\left( \frac{r}{2}\right)
;E\right) }\leq C\left( \eta ,p,q,\lambda ,E,A\right) r^{n\left( \frac{1}{q}-%
\frac{1}{p}\right) }\left\Vert \ w\right\Vert _{W^{2}{}^{,p}\left( Q\left(
r\right) ;E\right) }.  \tag{3.3}
\end{equation}

\bigskip \textbf{Proof. }\ Now we use $\left\Vert w\right\Vert _{2,\eta }$
to denote $\left\Vert w\right\Vert _{W^{2,\eta }\left( Q\left( r\right)
;E\right) }$ norm of $w\in W^{2,\eta }\left( Q\left( r\right) ;E\right) $
and we set $q\left( k\right) =\left( n+2\right) \left( n+2-k\right) ^{-1}$.
Now $\zeta $ is a nonnegative, compactly supported $C^{2,1}\left( Q\left(
r\right) \right) $ function with $\zeta \equiv 1$ in $Q\left( \frac{r}{2}%
\right) $ and 
\begin{equation*}
r\left\vert \nabla \zeta \right\vert +r^{2}\left\vert \nabla ^{2}\zeta
\right\vert +r\left\vert \zeta _{t}\right\vert \leq C\left( n\right) .
\end{equation*}

We then define 
\begin{equation*}
Q\left( w\right) =\oint_{Q\left( r\right) }w\left( x,t\right) dxdt,\text{ }%
Q\left( \nabla w\right) =\oint_{Q\left( r\right) }\nabla w\left( x,t\right)
dxdt,
\end{equation*}
\begin{equation*}
\bar{w}=w\left( x,t\right) -xQ\left( \nabla w\right) -Q\left( w\right)
\end{equation*}%
and set $w_{k}=\zeta ^{k}\bar{w}$. By resoning as in proof of Theorem 2.3 we
get 
\begin{equation*}
\left\Vert w_{k}\right\Vert _{2,q}\leq Cr^{-k}\left\Vert w\right\Vert _{2,1}
\end{equation*}%
as long as kon \c{s} 2 and the proof is completed by using Ho%
%TCIMACRO{\U{a8} }%
%BeginExpansion
\"{}
%EndExpansion
lder's inequality as before.

\bigskip From this estimate, we obtain the interior coercive estimate in
Morrey space.

\textbf{Theorem 3.2. }Let the assumptions (1)-(3) of Condition 1 be
satisfied. Assume $u\in W^{2,1,p}\left( Q\left( r\right) ;E\left( A\right)
,E\right) $ is a solution of the problem $\left( 1.4\right) $ and $f$ $\in
L^{p,\lambda }\left( Q\left( r\right) ;E\right) $ for some $\lambda \in
\left( 0,n+2\right) $. Then $\nabla ^{2}u\in L^{p,\lambda }\left( B\left( 
\frac{r}{2}\right) ;E\right) $, $\frac{\partial u}{\partial t}\in
L^{p,\lambda }\left( B\left( \frac{r}{2}\right) ;E\right) $ and there is a
constant $C=$ $C\left( \eta ,p,r,\lambda ,E,A\right) $ such that 
\begin{equation*}
\left\Vert \frac{\partial u}{\partial t}\right\Vert _{L^{p,\lambda }\left(
Q\left( \frac{r}{2}\right) ;E\right) }+\left\Vert \nabla ^{2}u\right\Vert
_{L^{p,\lambda }\left( Q\left( \frac{r}{2}\right) ;E\right) }+\left\Vert
Au\right\Vert _{L^{p,\lambda }\left( Q\left( \frac{r}{2}\right) ;E\right)
}\leq
\end{equation*}%
\begin{equation}
C\left( \left\Vert Q_{0}u\right\Vert _{L^{p,\lambda }\left( Q\left( r\right)
;E\right) }+\left\Vert \ \nabla ^{2}u\right\Vert _{L^{p}\left( Q\left(
r\right) ;E\right) }+\left\Vert Au\right\Vert _{L^{p}\left( Q\left( r\right)
;E\right) }\right) .  \tag{3.4}
\end{equation}

\textbf{Proof. }Let $\left( x_{0},t_{0}\right) \in Q\left( \frac{r}{2}%
\right) $ and let $\rho $ $>0$ be so small that $Q\left( x_{0},t_{0},\rho
\right) \subset Q\left( r\right) $. Let $q>p\left[ 1-\frac{\lambda }{n\left(
n+2\right) }\right] ^{-1}$, let $\upsilon \in W^{2,1,q}\left( Q\left(
x_{0},t_{0},\rho \right) ;E\right) $ be solution of the equation $\left(
1.4\right) $ in $Q\left( x_{0},t_{0},\rho \right) $ and $L_{1}\upsilon =0$
in $\partial Q\left( x_{0},t_{0},\rho \right) $. Finally, set $w=u-\upsilon $%
\ It follows from Theorem 3.1 that%
\begin{equation}
\left\Vert w\right\Vert _{W^{2,q}\left( Q\left( x_{0},t_{0},\rho \right)
;E\right) }\leq  \tag{3.5}
\end{equation}

\begin{equation*}
C\left( \eta ,p,q,\lambda ,E,A\right) \rho ^{\left( n+2\right) \left( \frac{1%
}{q}-\frac{1}{p}\right) }\left\Vert \ w\right\Vert _{W^{2,p}\left( Q\left(
x_{0},t_{0},\rho \right) ;E\right) }.
\end{equation*}%
The by reasoning as in the proof of Theorem 2.4 and by using Theorem A$_{3}$
we obtain the estimate $\left( 3.4\right) .$

\bigskip By reasoning as the above we obtain the corresponding boundary
estimate:

\textbf{Theorem 3.3. }Let the Condition 1 holds. Assume $u\in
W^{2,1,p}\left( Q^{+}\left( r\right) ;E\left( A\right) ,E\right) $ and $%
Lu\in $ $L^{p,\lambda }\left( Q\left( r\right) ;E\right) $. If also $%
L_{1}u=0 $ on $\partial Q\left( r\right) $, then $\nabla ^{2}u\in
L^{p,\lambda }\left( Q\left( \frac{r}{2}\right) ;E\right) $ and there is a
constant $C=$ $C\left( \eta ,p,r,\lambda ,E,A\right) $ such that 
\begin{equation*}
\left\Vert \frac{\partial u}{\partial t}\right\Vert _{L^{p,\lambda }\left(
Q^{+}\left( \frac{r}{2}\right) ;E\right) }+\left\Vert \nabla
^{2}u\right\Vert _{L^{p,\lambda }\left( Q^{+}\left( \frac{r}{2}\right)
;E\right) }+\left\Vert Au\right\Vert _{L^{p,\lambda }\left( Q^{+}\left( 
\frac{r}{2}\right) ;E\right) }\leq
\end{equation*}%
\begin{equation}
C\left( \left\Vert Lu\right\Vert _{L^{p,\lambda }\left( Q^{+}\left( r\right)
;E\right) }+\left\Vert \ \nabla ^{2}u\right\Vert _{L^{p}\left( Q^{+}\left(
r\right) ;E\right) }+\left\Vert Au\right\Vert _{L^{p}\left( Q^{+}\left(
r\right) ;E\right) }\right) .  \tag{3.6}
\end{equation}

\textbf{Theorem 3.4. \ }Let the Condition 1 holds for $\varphi >\frac{\pi }{2%
}$ and $0<\lambda <1$ and . Then for all $f\in \Phi $ problem $\left(
1.2\right) $ has a unique solution belonging to $\Phi ^{1,2}\left( A\right) $
and the following coercive estimate holds 
\begin{equation*}
\left\Vert \frac{\partial u}{\partial t}\right\Vert _{\Phi }+\left\Vert 
\frac{\partial ^{2}u}{\partial x^{2}}\right\Vert _{\Phi }+\left\Vert
Au\right\Vert _{\Phi }\leq C\left\Vert f\right\Vert _{\Phi }.
\end{equation*}

\bigskip \textbf{Proof. }\ Indeed, the prove Theorem 3.4 is obtained from
Theorems 3.1- 3.2 by localizasion and perturbation argument.

\bigskip

\begin{center}
\textbf{4. Application}

\textbf{4.1. Wentzell-Robin type mixed problem for elliptic equations}
\end{center}

\bigskip

\bigskip Consider the problem $\left( 1.3\right) -\left( 1.5\right) $ in
Morrey space $L^{\mathbf{p},\lambda }\left( \tilde{\Omega}\right) $. \ Let%
\begin{equation*}
\tilde{\Omega}=\mathbb{R}_{+}^{n}\times \left( 0,1\right) \text{, }X_{%
\mathbf{p}}=L^{\mathbf{p},\lambda }\left( \tilde{\Omega}\right) \text{, }X_{%
\mathbf{p}}^{2,\lambda }=W^{2,2,\lambda ,\mathbf{p}}\left( \left( \tilde{%
\Omega}\right) \right) \text{.}
\end{equation*}%
\begin{equation*}
\text{ }\mathbf{p=}\left( p_{1},p\right) ,\text{ }p_{1},p\in \left( 1,\infty
\right) .
\end{equation*}

\bigskip \textbf{Theorem 4.1. }Suppose\ $a$ is positive, $b$ is a
real-valued functions on $\tilde{\Omega}$. Moreover$,$ $a\left( x,.\right)
\in C\left[ 0,1\right] $ for all $x\in \mathbb{R}_{+}^{n}$, $a\left(
.,y\right) \in C_{b}\left( \mathbb{R}_{+}^{n}\right) $ for all $y\in \left[
0,1\right] $\ and $b\left( x,y\right) $ is a bounded function on $\tilde{%
\Omega}$ with 
\begin{equation*}
\exp \left( -\dint\limits_{\frac{1}{2}}^{x}b\left( x,y\right) a^{-1}\left(
x,y\right) dy\right) \in L_{1}\left( 0,1\right) \text{, for }x\in \mathbb{R}%
_{+}^{n};
\end{equation*}

(2) $a_{k}\in VMO\cap L^{\infty }\left( \mathbb{R}\right) ,$ $-a_{k}\left(
x_{k}\right) \in S\left( \varphi \right) ,$ $a_{k}\left( x_{k}\right) $ $%
\neq 0$ for a.e. $x_{k}\in \mathbb{R}$ and $\alpha _{i}\neq 0$, $i=0,1;$

(3) $\func{Re}\omega _{ki}\neq 0$ and $\frac{\lambda }{\omega _{ki}}\in
S\left( \varphi \right) $ for $\lambda \in S\left( \varphi \right) $, for
a.e. $x_{k}\in \mathbb{R}$, $0\leq \varphi <\pi ,$ $i=1,2$, $k=1,2,...,n$,
where $\omega _{ki}$ are roots of the equation $\left( 1.10\right) .$

\bigskip Then Then for all $f\in X_{\mathbf{p}}$, $\mu \in S\left( \varphi
\right) $ and for large enough $\left\vert \lambda \right\vert ,$ problem $%
\left( 1.3\right) -\left( 1.5\right) $ has a unique solution $u\in X_{%
\mathbf{p}}^{2,\lambda }.$ Moreover, the following coercive uniform estimate
holds%
\begin{equation}
\sum\limits_{k=1}^{n}\sum\limits_{i=0}^{2}\left\vert \mu \right\vert ^{1-%
\frac{i}{2}}\left\Vert \frac{\partial ^{i}u}{\partial x_{k}^{i}}\right\Vert
_{X_{\mathbf{p}}}+\left\Vert u\right\Vert _{X_{\mathbf{p}}}\leq C\left\Vert
\ f\right\Vert _{X_{\mathbf{p}}}.  \tag{4.1}
\end{equation}

\ \textbf{Proof.} Let $E=L^{p_{1}}\left( 0,1\right) $. It is known $\left[ 2%
\right] $\ that $L^{p_{1}}\left( 0,1\right) $ is an $UMD$ space for $%
p_{1}\in \left( 1,\infty \right) .$ Consider the operator $A$ defined by 
\begin{equation*}
D\left( A\right) =W^{2,p_{1}}\left( \tilde{\Omega};B_{j}u=0\right) ,\text{ }%
Au=a\frac{\partial ^{2}u}{\partial y^{2}}+b\frac{\partial u}{\partial y}.
\end{equation*}

Therefore, the problem $\left( 1.3\right) -\left( 1.5\right) $ can be
rewritten in the form of $\left( 1.1\right) $, where $u\left( x\right)
=u\left( x,.\right) ,$ $f\left( x\right) =f\left( x,.\right) $\ are
functions with values in $E=L^{p_{1}}\left( 0,1\right) .$ From $\left[ \text{%
10, 11}\right] $ we get that the operator $A$ generates analytic semigroup
in $L^{p_{1}}\left( 0,1\right) .$ Moreover, we obtain that the operator $A$
is $R$-positive in $L^{p_{1}}.$ Then from Theorem 2.1 we obtain the
assertion.

\begin{center}
\bigskip \textbf{4.2.} \textbf{The mixed problem for degenerate parabolic
equation}
\end{center}

\bigskip Consider the N-S problem $\left( 1.7\right) -\left( 1.9\right) $ in
Morrey space $L^{\mathbf{p},\lambda }\left( \tilde{\Omega}\right) $, where 
\begin{equation*}
\bar{\Omega}=\mathbb{R}_{+}^{n}\times \left( 0,b\right) \text{, }Y_{\mathbf{p%
}}=L^{\mathbf{p},\lambda }\left( \bar{\Omega}\right) \text{, }Y_{\mathbf{p}%
}^{2,1,\lambda }=W^{2,\left[ 2\right] ,1,\lambda ,\mathbf{p}}\left( \left( 
\bar{\Omega}\right) \right) \text{.}
\end{equation*}

The main aim of this section is to prove the following result:

\bigskip \textbf{Theorem 4.2. } Assume $b_{1}\left( .,y\right) \in C\left( 
\mathbb{R}_{+}^{n}\right) $ for all $y\in \left[ 0,b\right] ,$ $b_{1}\left(
x,.\right) \in C\left( \left[ 0,b\right] \right) $ for all $x\in \mathbb{R}%
_{+}^{n}$, $b_{1}\left( .,\right) \in C\left[ 0,1\right] $ and $b_{2}\left(
x,y\right) $ is a bounded function on $\bar{\Omega}$. Assume also $0\leq
\gamma <\frac{1}{2}$ and $\alpha _{10}\beta _{20}-\alpha _{20}\beta
_{10}\neq 0,$ $\alpha _{20}\beta _{11}+\alpha _{21}\beta _{10}-\alpha
_{10}\beta _{21}-\alpha _{11}\beta _{20}\neq 0,$ $\alpha _{11}\beta
_{21}+\alpha _{21}\beta _{11}\neq 0$, $\alpha _{11}\beta _{21}-\alpha
_{11}\alpha _{21}\neq 0,$ $\alpha _{11}\neq \beta _{11}$ for $\nu _{k}=1$
and $\left\vert \alpha _{k0}\right\vert +\left\vert \beta _{k0}\right\vert
>0,$ $\alpha _{10}\alpha _{20}+\beta _{10}\beta _{20}\neq 0$ for $\nu _{k}=0$%
. Moreover, suppose:

\ (1) $a_{k}\in VMO\cap L^{\infty }\left( \mathbb{R}\right) ,$ $-a_{k}\left(
x_{k}\right) \in S\left( \varphi \right) ,$ $a_{k}\left( x_{k}\right) $ $%
\neq 0$ for a.e. $x_{k}\in \mathbb{R}$ and $\alpha _{i}\neq 0$, $i=0,1;$

(2) $\func{Re}\omega _{ki}\neq 0$ and $\frac{\lambda }{\omega _{ki}}\in
S\left( \varphi \right) $ for $\lambda \in S\left( \varphi \right) $, for
a.e. $x_{k}\in \mathbb{R}$, $0\leq \varphi <\pi ,$ $i=1,2$, $k=1,2,...,n$,
where $\omega _{ki}$ are roots of the equation $\left( 1.10\right) .$

\bigskip Then Then for all $f\in Y_{\mathbf{p}}$, $\mu \in S\left( \varphi
\right) $ and for large enough $\left\vert \lambda \right\vert ,$ problem $%
\left( 1.7\right) -\left( 1.9\right) $ has a unique solution $u\in Y_{%
\mathbf{p}}^{2,1,\lambda }.$ Moreover, the following coercive uniform
estimate holds%
\begin{equation}
\sum\limits_{k=1}^{n}\sum\limits_{i=0}^{2}\left\vert \mu \right\vert ^{1-%
\frac{i}{2}}\left\Vert \frac{\partial ^{i}u}{\partial x_{k}^{i}}\right\Vert
_{Y_{\mathbf{p}}}+\left\Vert u\right\Vert _{Y_{\mathbf{p}}}\leq C\left\Vert
\ f\right\Vert _{Y_{\mathbf{p}}}.  \tag{4.2}
\end{equation}

\bigskip\ \textbf{Proof.} Let $E=L^{p_{1}}\left( 0,b\right) $. It is known $%
\left[ 3\right] $\ that $L^{p_{1}}\left( 0,b\right) $ is an $UMD$ space for $%
p_{1}\in \left( 1,\infty \right) .$ Consider the operator $A$ defined by 
\begin{equation*}
D\left( A\right) =W^{\left[ 2\right] ,p_{1}}\left( \Omega ;L_{k}u=0\right) ,%
\text{ }Au=b_{1}\frac{\partial ^{\left[ 2\right] }u}{\partial y^{2}}+b_{2}%
\frac{\partial ^{\left[ 1\right] }u}{\partial y}.
\end{equation*}

Therefore, the problem $\left( 1.7\right) -\left( 1.9\right) $ can be
rewritten in the form of $\left( 1.4\right) $, where $u\left( x\right)
=u\left( x,.\right) ,$ $f\left( x\right) =f\left( x,.\right) $\ are
functions with values in $E=L^{p_{1}}\left( 0,b\right) .$ From $\left[ 26%
\right] $ we get that the operator $A$ generates analytic semigroup in $%
L^{p_{1}}\left( 0,b\right) .$ Moreover, we obtain that the operator $A$ is $%
R $-positive in $L^{p_{1}}$. Then from Theorem 3.4 we obtain the assertion.

\bigskip

%\begin{center}
%\bigskip
%
%4. \textbf{Applications}
%\end{center}

\section*{Acknowledgements}

The first author is partially supported by P.R.I.N. 2019 and the ``RUDN
University Program 5-100''. \bigskip

\begin{center}
\textbf{References}
\end{center}

\begin{enumerate}
\item Amann H., Linear and quasi-linear equations,1, Birkhauser, 1995.

\item Bourgain J., Some remarks on Banach spaces in which martingale
differences are unconditional, Arkiv Math. 21 (1983), 163-168.

\item Burkholder D. L., A geometrical conditions that implies the existence
certain singular integral of Banach space-valued functions, Proc. conf.
Harmonic analysis in honor of Antonu Zigmund, Chicago, 1981,Wads Worth,
Belmont, (1983), 270-286.

\item Besov, O. V., P. Ilin, V. P., Nikolskii, S. M., Integral
representations of functions and embedding theorems, Nauka, Moscow, 1975.

\item Denk R., Hieber M., Pr\"{u}ss J., $R$-boundedness, Fourier multipliers
and problems of elliptic and parabolic type, Mem. Amer. Math. Soc. 166
(2003), n.788.

\item Favini A., Shakhmurov V., Yakubov Y., Regular boundary value problems
for complete second order elliptic differential-operator equations in UMD
Banach spaces, Semigroup Form, v. 79 (1), 2009.

\item Chiarenza F., Frasca M., and Longo P., Interier $W^{2,p}$estimates for
non divergence elliptic equations with discontinuous coefficients, Ricerche
Mat. 40, 1991, 149-168.

\item Chiarenza F., Frasca M., and Longo P., $W^{2,p}$-solvability of the
Dirichlet problem for nondivergence elliptic equations with VMO
coefficients, Transactions of the American Mathematical Socity, 336(2),1993,
841-853.

\item P. Cannarsa, Second order non variational parabolic systems, Boll. Un.
Mat. Ital. 18-C (1981), 291--315.

\item A. Favini, G. R. Goldstein, Jerome A. Goldstein and Silvia Romanelli,
Degenerate Second Order Differential Operators Generating Analytic
Semigroups in $L_{p}$ and $W^{1,p}$, Math. Nachr. 238 (2002), 78 -- 102.

\item V. Keyantuo, M. Warma, The wave equation with Wentzell-Robin boundary
conditions on Lp-spaces, J. Differential Equations 229 (2006)
680--697.Gorbachuk V. I. and Gorbachuk M. L., Boundary value problems \ for
differential-operator equations, Naukova Dumka, Kiev, 1984

\item L.C. Evans, R.P. Gariepy, Measure Theory and Fine Properties of
Functions, CRC Press, Boca Raton, 1992.

\item G. DiFazio, D.K. Palagachev, Oblique derivative problem for elliptic
equations in non divergence form with VMO coefficients, Comment. Math. Univ.
37 (1996) 537--556.

\item G. Di Fazio, D.K. Palagachev, M.A. Ragusa, Global Morrey regularity of
strong solutions to Dirichlet problem for elliptic equations with
discontinuous coefficients, J. Funct. Anal. 166 (1999)176--196.

\item O. A. Ladyzenskaja, V.A. Solonnikov, N.N. Ural0ceva, Linear and
Quasi-linear Equations of Parabolic Type, American Mathematical Society,
Providence, RI, 1968.

\item G. M. Liebermann, \ A mostly elementary proof of Morrey space
estimates for elliptic and parabolic equations with VMO coefficients, J.
Func. Analysis, 201 (2003), 457-479.

\item G. M. Lieberman, The natural generalization of the natural conditions
of Ladyzhenskaya and Uraltseva for elliptic equations, Comm. Partial
Differential Equations 16 (1991) 311--361.

\item A. Maugeri, D.K. Palagachev, Boundary value problems with an oblique
derivative for uniformly elliptic operators with discontinuous coefficients,
Forum Math. 10 (1998) 393--405.

\item M. A. Ragusa, Commutators of fractional integral operators on
vanishing-Morrey spaces, J. Global Optim. 40(1--3) (2008) 361--368.

\item M. A. Ragusa, Homogeneous Herz spaces and regularity results,
Nonlinear Anal. 71(12) (2009) e1909--e1914.

\item M. A. Ragusa, Embeddings for Morrey--Lorentz spaces, J. Optim. Theory
Appl. 154 (2) (2012) 491--499.

\item M. A. Ragusa, A. Scapellato, Mixed Morrey spaces and their
applications to partial differential equations, Nonlinear Analysis 151
(2017) 51--65.

\item M. A. Ragusa, V. B. Shakhmurov, \ Embedding of vector-valued Morrey
spaces and separable differential operators, Bulletin of Mathematical
Sciences 9(2) (2019) 1950013, DOI: 10.1142/S1664360719500139.

\item Shakhmurov V. B., Embedding theorems and maximal regular differential
operator equations in Banach-valued function spaces, Journal of Inequalities
and Applications, ( 4 ) 2005, 605-620.

\item Shakhmurov V. B., Coercive boundary value problems for regular
degenerate differential-operator equations, J. Math. Anal. Appl., 292 ( 2),
(2004), 605-620.

\item Shakhmurov V. B., Separable anisotropic differential operators and
applications, J. Math. Anal. Appl. 2007, 327(2), 1182-1201.

\item Shakhmurov V. B., Abstract elliptic equations with VMO coefficients in
half plane, Mediterranen Journal of Mathematics, 25(2015), 1-21.

\item Weis L., Operator-valued Fourier multiplier theorems and maximal $%
L_{p} $ regularity, Math. Ann. 319, (2001), 735-75.

\item Triebel H., Interpolation theory. Function spaces. Differential
operators., North-Holland, Amsterdam, 1978.

\item Yakubov S. and Yakubov Ya., Differential-operator Equations. Ordinary
and Partial \ Differential Equations, Chapman and Hall /CRC, Boca Raton,
2000.

\item W. P. Ziemer, Weakly Differentiable Functions, Springer, New York,
1989.
\end{enumerate}

\bigskip

\end{document}